\documentclass[11pt, letterpaper, reqno]{amsart}
\usepackage{t1enc}
\usepackage[latin1]{inputenc}
\usepackage[english]{babel}
\usepackage{amsmath,amsthm}
\usepackage{amsfonts}
\usepackage{latexsym}
\usepackage[dvips]{graphicx}
\usepackage{float}
\usepackage[natural]{xcolor}
\usepackage{algpseudocode}
\usepackage{amssymb}
\usepackage{enumerate}
\usepackage{enumitem}
\usepackage{multirow}
\usepackage{xcolor}
\usepackage[colorlinks,linkcolor=blue]{hyperref}
\usepackage{lineno}
\usepackage{setspace}
\usepackage{fullpage}
\usepackage{bm}
\usepackage{lineno}
\usepackage{indentfirst}
\usepackage[letterpaper, margin=1in]{geometry} 
\usepackage{setspace}
\usepackage[normalem]{ulem}
\usepackage{listings}
\usepackage[ruled, linesnumbered, vlined]{algorithm2e}
\SetKwInOut{Input}{Input}
\SetKwInOut{Output}{Output}
\SetKw{KwRet}{return}
\theoremstyle{plain}

\numberwithin{equation}{section}
\newtheorem{thm}{Theorem}[section]
\newtheorem{theorem}[thm]{Theorem}
\newtheorem{conjecture}[thm]{Conjecture}
\newtheorem{lemma}[thm]{Lemma}
\newtheorem{corollary}[thm]{Corollary}
\newtheorem{proposition}[thm]{Proposition}

\theoremstyle{definition}

\newtheorem{question}[thm]{Question}

\newtheorem{remark}[thm]{Remark}

\newtheorem{definition}[thm]{Definition}
\newtheorem{claim}[thm]{Claim}
\newtheorem{fact}[thm]{Fact}

\newtheorem{example}[thm]{Example}

\newtheorem{defn-thm}[thm]{Definition-Theorem}

\newcommand{\btheorem}{\begin{theorem}}
	\newcommand{\etheorem}{\end{theorem}}
\newcommand{\bconjecture}{\begin{conjecture}}
	\newcommand{\econjecture}{\end{conjecture}}
\newcommand{\bproposition}{\begin{proposition}}
	\newcommand{\eproposition}{\end{proposition}}
\newcommand{\bdefinition}{\begin{definition}}
	\newcommand{\edefinition}{\end{definition}}
\newcommand{\bcorollary}{\begin{corollary}}
	\newcommand{\ecorollary}{\end{corollary}}
\newcommand{\bproof}{\begin{proof}}
	\newcommand{\eproof}{\end{proof}}
\newcommand{\bclaim}{\begin{claim}}
	\newcommand{\eclaim}{\end{claim}}
\newcommand{\bquestion}{\begin{question}}
	\newcommand{\equestion}{\end{question}}
\newcommand{\bfact}{\begin{fact}}
	\newcommand{\efact}{\end{fact}}
\newcommand{\bremark}{\begin{remark}}
	\newcommand{\eremark}{\end{remark}}
\newcommand{\eexample}{\end{example}}
\newcommand{\bexample}{\begin{example}}

\newcommand{\elemma}{\end{lemma}}
\newcommand{\blemma}{\begin{lemma}}

\newcommand{\beq}{\begin{equation}}
	\newcommand{\eeq}{\end{equation}}

\newcommand{\mathify}[1]{\ifmmode{#1}\else\mbox{$#1$}\fi}
\newcommand{\bigO}O

\newcommand\card[1]{\left| #1 \right|}
\newcommand\cardi[1]{| #1 |}

\renewcommand{\vec}[1]{{\mathbf #1}}

\newcommand{\remove}[1]{{}}

\begin{document}

\title{Perfect Matchings in Random Sparsifications of Dense Hypergraphs}
\author{Jie Han}
\author{Jingwen Zhao}
\address{School of Mathematics and Statistics and Center of Applied Mathematics, Beijing Institute of Technology, Beijing, China, 100081. 
Email: {\tt \{han.jie|jingwen.zhao\}@bit.edu.cn}.}

\begin{abstract}
Given \(1\le\ell <k \) and \(\delta\geq 0\), let \(\mathbf{PM}(k,\ell,\delta)\) be the decision problem for the existence of perfect matchings in \(n\)-vertex \(k\)-uniform hypergraphs with minimum \(\ell\)-degree at least \(\delta\binom{n-\ell}{k-\ell}\). For \(k\geq 3\), \(\mathbf{PM}(k,\ell,0)\) was one of the first NP-complete problems identified by Karp. Keevash, Knox and Mycroft conjectured that \(\mathbf{PM}(k,\ell,\delta)\) is in P for every \(\delta>1-(1-1/k)^{k-\ell}\) and this was recently verified by the work of Gan--Han, together with a very recent work of Fu et al.

In this paper we study the existence of perfect matchings in the random $p$-sparsification of such $k$-uniform hypergraphs, that is, for $p=p(n)\in [0,1]$, each edge is selected independently with probability \(p\). 
Building on the structural theory of Gan and Han, we show that the corresponding dense perfect matching results are robust under random sparsification. 
As consequences, we obtain deterministic polynomial-time algorithms that asymptotically almost surely solve the associated decision problems, as well as lower bounds on the number of perfect matchings in such hypergraphs -- interestingly, such hypergraphs either have no perfect matching, or have $(\Omega(n))^{(1-1/k)n}$ perfect matchings. 
Moreover, we also establish analogous results for the \(F\)-factor problem in graphs.

Our proofs combine a partial exposure algorithm, the lattice-based absorption method, and a random redistribution method of Kelly, M\"uyesser and Pokrovskiy, via the framework of spread distributions. 
A key new ingredient is a lattice-preparation step that separates the contributions of the two classes of robust index vectors arising in the Gan--Han structural theory. 
Together with the random redistribution method, this allows us to establish the desired spread property in the family of perfect matchings.
\end{abstract}

\maketitle

\section{Introduction}

Matchings are fundamental objects in Graph Theory and have broad applications in other branches of Science and a variety of practical problems
(e.g.\ the assignment of graduating medical students to their first hospital appointments\footnote{In 2012, the Nobel Memorial Prize in Economics was awarded to Shapley and Roth ``for the theory of stable allocations and the practice of market design."}).
Applications of matchings in hypergraphs 
include the `Santa Claus' allocation problem~\cite{santa_claus};
they also offer a universal framework for many 
important combinatorial problems, e.g.\ the Existence Conjecture 
for designs (see \cite{Keevash_design, GKLO}) and
Ryser's conjecture \cite{ryser} on transversals in Latin squares.

Given $k\ge 2$, a $k$-uniform hypergraph (or simply a \emph{$k$-graph}) is a hypergraph $H = (V(H), E(H))$ in which every edge $e \in E(H)$ has size exactly $k$.
A \emph{matching} $M$ in $H$ is a set of vertex-disjoint edges. 
A \emph{perfect matching} in $H$ is a matching that covers all the vertices of $H$.
    
Determining the existence of perfect matchings in $k$-graphs is a central problem in combinatorics.  
Several well-known conjectures and problems, such as Ryser's conjecture and the existence of certain combinatorial designs, can be reduced to this type of problems. 
When $k=2$, Tutte \cite{tutte1947factorization} provided a complete characterization of graphs that contain a perfect matching in 1947, and Edmonds \cite{edmonds1965paths} gave a polynomial-time algorithm to determine whether a given graph contains one. 
However, for $k \ge 3$, the problem becomes \(\mathsf{NP}\)-complete, as shown by Karp \cite{karp2009reducibility}. 
Naturally, we aim to find sufficient conditions that guarantee the existence of perfect matchings. 
Among these, the most extensively studied is the minimum degree condition, commonly referred to as the Dirac-type problem.

As an important generalization of perfect matchings, the concept of \emph{$F$-factors} has received considerable attention in graph and hypergraph theory. 
Given two $k$-graphs $F$ and $H$, an \emph{$F$-packing} of $H$ is a collection of vertex-disjoint copies of $F$ in $H$. 
An \emph{$F$-factor}, or \emph{perfect \(F\)-packing}, is an $F$-packing that covers all vertices of $H$.
Hell and Kirkpatrick \cite{kirkpatrick1983complexity} showed that determining whether a graph $G$ contains a perfect $F$-packing is \(\mathsf{NP}\)-complete when $F$ has a component with at least three vertices.


Besides studying the containment of an \(F\)-factor, it is also interesting to ask how stable this property is. 
In other words, does the \(F\)-factor still exist if the graph is slightly changed? 
This leads to the study of the robustness of graph properties. 
Several measures of robustness have been proposed. 
For example, one can measure the robustness of Dirac graphs (i.e., graphs \(G\) with minimum degree at least \(\card{G}/2\)) with respect to perfect matchings by computing the number of perfect matchings that a Dirac graph must contain. 
Indeed, confirming a conjecture of S\'{a}rk\"{o}zy,  Selkow, and Szemer\'{e}di\cite{sarkozy2003number}, Cuckler and Kahn \cite{cuckler2009hamiltonian} proved that every Dirac graph contains at least \(n!/{(2+o(1))}^n\) perfect matchings.

Another measure is the so-called resilience. 
Roughly speaking, for monotone increasing graph properties, the resilience quantifies the robustness in terms of the number of edges one must delete from \(H\), locally or globally, in order to destroy the property \(\mathcal{P}\). 
For example, Krivelevich, Lee and Sudakov \cite{krivelevich2014robust} proved that for any Dirac graph \(H\), the random subgraph obtained by keeping each edge of \(H\) independently with probability \(p\) where \( p\ge C \log n/n\) (\(C\) is an absolute constant) \emph{asymptotically almost surely} (abbreviated a.a.s.) contains a perfect matching.
    
In this work, we focus on the robustness of degree conditions that guarantee perfect matchings or $F$-packings in hypergraphs. 
While such conditions have traditionally been studied from a structural or algorithmic perspective (particularly in the context of making the decision problem polynomial-time solvable), we show that many of them also exhibit a form of robustness.
That is, they continue to guarantee the desired structures even under random perturbations, such as when passing to a random subgraph.
In particular, we identify the two settings in which the same minimum degree thresholds that make the problem tractable also lead to robust containment in random subgraphs.

 \subsection{Perfect matchings in hypergraphs}

    Let \(k, \ell \in \mathbb{N}\) where \(\ell \le k-1\). 
    Given a $k$-graph $H$ and a subset $S\subseteq V(H)$ of size $\ell$, the \emph{degree} of $S$ in $H$, denoted by $\deg_H(S)$, is the number of edges in $H$ that contain $S$. 
    The minimum $\ell$-degree $\delta_\ell(H)$ is the minimum degree over all $\ell$-sets of $V(H)$. 
    
    One of the central questions in the field is to determine minimum degree conditions that guarantee the existence of a perfect matching. 
    When the degree condition is sufficiently strong (e.g., $1 - 1/k$, as follows from a result in \cite{han2009perfect}), the existence of a perfect matching is guaranteed, and the corresponding decision problem becomes trivial.
    From an algorithmic perspective, it is therefore natural to study for which $\ell$-degree thresholds the perfect matching decision problem remains computationally tractable. 
    To formalize this, let \(\mathbf{PM}(k,\ell,\delta)\) denote the decision problem of determining whether a \(k\)-graph \(H\) on \(n\) vertices with \(\delta_\ell(H) \ge \delta \binom{n-\ell}{k-\ell}\) contains a perfect matching. 
    Szyma\'{n}ska \cite{szymanska2013complexity} proved that \(\mathbf{PM}(k,\ell,\delta)\) is \(\mathsf{NP}\)-complete whenever \(\delta < 1 - (1 - 1/k)^{k-\ell}\), by a polynomial-time reduction from \(\mathbf{PM}(k,\ell,0)\).
    This motivated Keevash, Knox and Mycroft \cite{keevash2013polynomial} to propose the following conjecture, which they verified in the case \(\ell = k-1\). 
    \begin{conjecture}[Keevash, Knox and Mycroft \cite{keevash2013polynomial}]\label{conj:keevash}
        \(\mathbf{PM}\left(k,\ell,\delta\right)\) is in \(\mathsf{P}\) for \(\delta>1-(1-1/k)^{k-\ell}\).
    \end{conjecture}
    Subsequently, Han and Treglown \cite{han2020complexity} showed that the conjecture holds for \(0.5k\le \ell \le (1+\ln(2/3))k\approx  0.59 k\). 
    Recently, Gan and Han~\cite{gan2025keevash} made significant progress by reducing the problem to minimum degree thresholds for perfect fractional matchings.  
    Given a \(k\)-graph \(H = (V, E)\), a \emph{fractional matching} in \(H\) is a function \(w :E\to [0,1]\) such that for each \(v\in V(H)\) we have that \(\sum_{e\ni v}w(e)\le 1 \). Then \(\sum_{e\in E(H)}w(e)\) is the \emph{size} of \(w\). If the size of \(w\) in \(H\) is \(n/k\) then we say that \(w\) is a \emph{perfect fractional matching}. Given \(k,\ell\in\mathbb{N}\) such that \(\ell\le k-1\), define \(c^*_{k,\ell}\) to be the smallest number \(c\) such that every \(k\)-graph \(H\) on \(n\) vertices with \(\delta_\ell(H)\ge (c+o(1))\binom{n-\ell}{k-\ell}\) contains a perfect fractional matching. 

    \begin{theorem}\cite[Theorem 1.3]{gan2025keevash}\label{thm:gan_han_algorithm}
        Suppose \( k, \ell \in \mathbb{N} \) with \( 1 \leq \ell \leq k-1 \). Then for any \( \delta \in (c_{k,\ell}^*, 1] \), \( \mathbf{PM}(k, \ell, \delta) \) is in \( \mathsf{P} \). That is, for any \( \delta \in (c_{k,\ell}^*, 1] \), there exists a constant \( c = c(k) \) such that there is an algorithm with running time \( O(n^c) \) which given any \( n \)-vertex \( k \)-graph \( H \) with \( \delta_\ell(H) \geq \delta \binom{n-\ell}{k-\ell} \), determines whether \( H \) contains a perfect matching.
    \end{theorem}

    They in fact proved a stronger result, providing a polynomial-time algorithm that outputs either a perfect matching or a certificate that none exists. 
    For the parameter \(c_{k,\ell}^*\), Alon, Frankl, Huang, R\"{o}dl, Ruci\'{n}ski, and Sudakov \cite{alon2012large} made the following conjecture, and verified the case \(k-\ell \le 4\). 
    \begin{conjecture}\cite[Conjecture 1.1]{alon2012large}
    For all $k,\ell \in \mathbb{N}$ with $\ell \le k-1$, \(c^*_{k,\ell} = 1 - (1 - 1/k)^{k-\ell}\).
    \end{conjecture}
    The conjecture is known to hold in several ranges of parameters. 
    K\"{u}hn, Osthus and Townsend \cite{kuhn2014fractional} proved it for \(\ell\ge k/2\), and Han \cite{han2016perfect} established the case \(\ell=(k-1)/2\). 
    More recently, Frankl and Kupavskii \cite{frankl2022erdHos} verified the conjecture for \(\ell\ge 0.4k\). 
    A very recent resolution of a conjecture of Feige \cite{fu2026sharp} together with a result of Ferber and Jain, implies the full validity of Conjecture 1.3.
    Combining with Theorem~\ref{thm:gan_han_algorithm}, it follows that Conjecture~\ref{conj:keevash} holds.
    
    In this paper, we strengthen this line of work by establishing a robust version of Theorem \ref{thm:gan_han_algorithm}. 
    The conference version of this work established the following theorem in the special case $\ell=k-1$, while the present paper extends it to all $1\le \ell \le k-1$. 
    For any hypergraph \(H\) and \(0\le p\le 1\), let \(H_p\) be a spanning random subhypergraph of \(H\) obtained by choosing each edge \(e\in H\) with probability \(p\) independently at random.

    \begin{theorem}[Main result]\label{thm:algorithm_PM}
       Suppose \( k, \ell \in \mathbb{N} \) with \( 1 \leq \ell \leq k-1 \). 
       Then for any \(\delta\in (c_{k,\ell}^*,1] \), there exist \(c=c(k)\) and \(C=C(k,\delta)>0\) such that there is an algorithm with running time \(O(n^c)\) which, given any \(n\)-vertex \(k\)-graph \(H\) with \(\delta_\ell(H)\ge \delta \binom{n-\ell}{k-\ell}\), has the following properties for \(p\ge C\log n / n^{k-1} \):
        \begin{itemize}
        \item If the algorithm accepts, then a.a.s.~$H_p$ contains a perfect matching;
        \item If the algorithm rejects, then $H_p$ does not contain a perfect matching.
        \end{itemize}
    \end{theorem}

This result establishes a stronger, probabilistic guarantee concerning rather sparse random subgraphs of the original $k$-graph. 
Moreover, it appears to be the first robust (random-sparsification) version in the polynomial-time tractable regime of the perfect matching decision problem. 
The assumptions and conclusions are essentially best possible in multiple senses: 
\begin{itemize}
    \item for $p=o(\log n/n^{k-1})$, $H_p$ a.a.s.~contains no perfect matching even for complete $k$-graph $H$;
    \item if $\delta_{\ell}(H) \le (c_{k,\ell}^* - \gamma)\binom{n-\ell}{k-\ell}$ then as mentioned above, Szyma\'nska's result implies that the decision problem for perfect matching in $H$ is $\mathsf{NP}$-complete;
    \item the one-sided probabilistic guarantee in the conclusion is also the best one can hope; indeed, a (deterministic) algorithm deciding the existence of perfect matchings even for the random sparsification of the complete $k$-graph would not exist unless \(\mathsf{P}=\mathsf{NP}\).
\end{itemize}

For the decision problem \(\mathbf{PM}\left(k,\ell,\delta\right)\), significant challenges arise for the general case $\ell < k-1$, because of the (possible) existence of a small set of exceptional vertices.
Building on previous lattice absorption theory, 
Gan and Han \cite{gan2025keevash} resolved the decision problem by introducing two rather different types of robust vectors to classify these vertices.
In our setting, we have to deal with these two types of vectors very differently.
This motivates us to introduce a lattice-preparation step, implemented by a preliminary matching with prescribed index vectors, which reduces the two-class lattice structure to the one needed for the later random redistribution argument while keeping the random choices sufficiently controlled.

A related, though weaker, robust version for \(F\)-factors will be presented in Section~\ref{sec:F-factor}.
Now we turn to perfect $F$-packings in graphs and state our second main result, Theorem \ref{thm:main_thm_factor}, after some introduction and preparation.

\subsection{Perfect packings in graphs}

    The Hajnal-Szemer\'{e}di theorem states that every \(n\)-vertex graph with minimum degree at least \((1-1/r)n\) and \(r\mid n\) contains a \(K_r\)-factor, and this bound is tight. More generally, Alon and Yuster~\cite{alon1996h} proved the following theorem. Given a graph \(F\), let \(v_F\) denote its number of vertices, \(e_F\) its number of edges, and \(\chi (F)\) its chromatic number. Let \(M(n)\) be the time needed to multiply two $n$ by $n$ matrices with \(0,1\) entries. Determining $M(n)$ is a challenging problem in theoretical computer science, and the best known bound of $M(n) = O(n^{2.3728596})$ was obtained recently by Alman and Williams~\cite{alman2024refined}.

    \begin{theorem}\cite[Theorem 1.1]{alon1996h}\label{alon-chromatic number}
        For every \(\gamma>0\) and each graph \(F\) there exists an integer \(n_0=n_0(\gamma,F)\) such that every graph \(G\) whose order \(n>n_0\) is divisible by \(v_F\) and whose minimum degree is at least \(\left(1-1/\chi(F)+\gamma\right)n\) contains an \(F\)-factor. Moreover, there is an algorithm which finds such \(F\)-factor in time \(O(M(n))\).
    \end{theorem}

    Alon and Yuster also provided examples showing that the error term \(\gamma n\) cannot be completely omitted for some graphs, but they conjectured that it could be replaced by a constant depending only on \(F\). 
    This conjecture was later proved by Koml\'{o}s, S\'{a}rk\"{o}zy and Szemer\'{e}di~\cite{komlos2001proof}.

    On the other hand, there are graphs for which the bound on the minimum degree can be significantly improved by replacing the chromatic number with a refined parameter known as the \emph{critical chromatic number}~\cite{cooley2007perfect}.
    The \emph{critical chromatic number} \(\chi_{cr}(F)\) of a graph \(F\) is defined as \(\left(\chi(F)-1\right)v_F/\left(v_F-\sigma(F)\right)\), where \(\sigma(F)\) denotes the  minimum size of the smallest colour class in a colouring of \(F\) with \(\chi(F)\) colours. 
    Note that \(\chi(F)-1<\chi_{cr}(F)\le \chi(F)\) and the equality holds if and only if every \(\chi(F)\)-colouring of \(F\) has equal colour class sizes.  
    If $\chi_{cr}(F)=\chi(F)$, then we call $F$ \emph{balanced}; otherwise we call it \emph{unbalanced}.
    Koml\'{o}s~\cite{komlos2000tiling} initiated the study of $F$-packings at the threshold determined by the critical chromatic number. 
    This was later strengthened by Shokoufandeh and Zhao~\cite{shokoufandeh2003proof}, who showed that the leftover can be bounded by a constant depending only on $F$.

    \begin{theorem}\cite[Theorem 1.8]{shokoufandeh2003proof}\label{critical chromatic number-alomst pefect F-factor}
        For any $F$ there is an $n_0=n_0(F)$ so that if $G$ is a graph on $n\ge n_0$ vertices and minimum degree at least $\left(1-1/\chi_{cr}(F)\right)n$, then $G$ contains an $F$-packing that covers all but at most $5{v_F}^2$ vertices.
    \end{theorem}

    Given the two types of minimum degree thresholds discussed above, it is natural to consider the corresponding decision problem. 
    Let \(\mathbf{Pack}(F, \delta)\) denote the problem of determining whether a graph \(G\) with minimum degree \(\delta(G)\ge\delta \card{G}\) contains an \(F\)-factor. 
    As mentioned earlier, when \(F\) contains a component with at least three vertices, results of Hell and Kirkpatrick~\cite{kirkpatrick1983complexity} show that \(\mathbf{Pack}(F,0)\) is \(\mathsf{NP}\)-complete. 
    More generally, K\"{u}hn and Osthus~\cite{KuehnOsthus2006} proved that \(\mathbf{Pack}(F,\delta)\) is \(\mathsf{NP}\)-complete for any \(\delta \in [0,1-1/\chi_{cr}(F))\), provided that \(F\) is either a clique of size at least \(3\) or a complete \(k\)-partite graph where \(k\ge 2\) and the second smallest vertex class has size at least \(2\). 
    In contrast, Theorem~\ref{alon-chromatic number} implies that \(\mathbf{Pack}(F,\delta)\) is in \(\mathsf{P}\) for all \(\delta \in (1-1/\chi(F), 1]\). 
    Recently, Han and Treglown~\cite{han2020complexity} provided a polynomial-time algorithm showing that \(\mathbf{Pack}(F,\delta)\) is in \(\mathsf{P}\) whenever \(\delta \in (1 - 1/\chi_{cr}(F),1]\), answering a question of Yuster negatively.

    \begin{theorem}\cite[Theorem 1.11]{han2020complexity}\label{algorithm-factor}
        For any \(r\)-vertex \(k\)-chromatic graph \(F\) and \(\delta\in (1-1/\chi_{cr}(F), 1]\), \(\mathbf{Pack}(F, \delta)\) is in \(\mathsf{P}\). That is, for every \(n\)-vertex graph \(G\) with minimum degree at least \(\delta n\), there is an algorithm with running time \(O\left(n^{\max \{2^{r^{k-1}-1}r+1, r(2r-1)^r\}}\right)\), which determines whether \(G\) contains a perfect \(F\)-packing.
    \end{theorem}

    Theorem~\ref{algorithm-factor} can efficiently detect the divisibility obstructions when the minimum degree condition guarantees an \(F\)-packing that covers all but a constant number of vertices (as established in Theorem~\ref{critical chromatic number-alomst pefect F-factor}). 
    Similar to the perfect matching problem, it is natural to seek a robust version of Theorem~\ref{algorithm-factor} similar to Theorem \ref{thm:algorithm_PM}. 
    Unfortunately, we are only able to prove a partial result which we state below after some preparations.

\subsection{\texorpdfstring{Partition, lattices and new result on $F$-packings}{Partition, lattices and new result on F-packings}}\label{sec:F-factor}

The proofs in \cite{gan2025keevash} and \cite{han2020complexity} used a lattice-based absorbing method, which combines the absorbing technique with the concept of divisibility barriers, introduced in \cite{keevash2015geometric}. 
    It was shown in \cite{keevash2015geometric} that a $k$-graph $H$ either contains a perfect matching or is structurally close to a family of lattice-based constructions known as divisibility barriers. 
    To describe divisibility barriers in general, we begin with the following definition. 
    In this paper, every partition has an implicit order on its parts. 

    \begin{definition}
    [Partition, index vector and lattice]\label{def_lattice}
        Let \(H=\left(V,E\right)\) be an \(n\)-vertex \(k\)-graph, and let \(\mathcal{P}=\{V_0,V_1,\dots,V_d\}\) be a partition of \(V\), where \(V_0\) is allowed to be empty.  
        Let \(F\) be an \(r\)-vertex \(k\)-graph. 
        The \emph{index vector} \(\vec{i}_{\mathcal{P}}(S)\in \mathbb{Z}^d\) of a subset \(S\subseteq V\) with respect to \(\mathcal{P}\) is defined by \(\vec{i}_{\mathcal{P}}(S)|_i=\card{S\cap V_i}\) for each \(i\in[d]\), where \(\vec{v}|_i\) denotes the \(i\)th coordinate of \(\vec{v}\). 
        Note that for index vectors, the part \(V_0\) is not considered.
        Then for any \(\mu >0\), 
        \begin{enumerate}[label=(\arabic*)]
            \item \(I_{\mathcal{P},F}^{\mu}(H)\) denotes the set of \(\vec{i}\in \mathbb{Z}^d\) such that \(H\) contains at least \(\mu n^r\) copies of \(F\) with index vector \(\vec{i}\);
            \item \(L_{\mathcal{P},F}^{\mu}(H)\) denotes the lattice (that is, the additive subgroup) in \(\mathbb{Z}^d\) generated by \(I_{\mathcal{P},F}^{\mu}(H)\).
        \end{enumerate}
    \end{definition}
    
    We call a copy of \(F\) \emph{\(\mu\)-robust} if its index vector $\vec{i}\in I_{\mathcal{P},F}^{\mu}(H)$.
    For each vertex $ v\in V(H) $, let $ F^{\mu}_{\mathcal{P}}(v) $ be the collection of \(r\)-sets that contain $v$ and span at least one \(\mu\)-robust copy of \(F\). 
    When \(F\) is a single edge, we abbreviate these as \(I_{\mathcal{P}}^{\mu}(H)\), \(L_{\mathcal{P}}^{\mu}(H)\), and \( E^{\mu}_{\mathcal{P}}(v) \), respectively.

    To apply the absorbing method, we first require a suitable partition, which is provided by the following lemma. Before stating it, we give some basic definitions.
    We use the reachability arguments introduced by Lo and Markstr\"{o}m \cite{lo2015f}. 
    We say that two vertices \(u\) and \(v\) in \(V(H)\) are \((F,\beta,t)\)-\emph{reachable} in \(H\) if there are at least \(\beta n^{tr - 1}\) \((tr-1)\)-sets \(S\) such that both \( H[S\cup\{u\}] \) and \( H[S\cup\{v\}] \) have perfect \(F\)-packings. 
    We refer to such a set \(S\) as a \emph{reachable} \((tr-1)\)\emph{-set for \( u \) and \( v \)}. 
    A set \(U \subseteq V(H)\) is \((F,\beta,t)\)-\emph{closed} if every pair of vertices in \(U\) is \((F,\beta,t)\)-reachable. 
    We say that a partition \(\mathcal{P}=\{V_1,\ldots,V_d\}\) of \(V(H)\) is \((F,\beta,t,c)\)-\emph{good} if each part \(V_i\) is \((F,\beta,t)\)-closed and \(|V_i|\ge cn\) for all \(i\in[d]\). 
    When \(F\) is a single edge, we simply write \((\beta,t)\)-reachable, \((\beta,t)\)-closed, and \((\beta,t,c)\)-good.
    Throughout the paper, we write \(\alpha\ll \beta\ll\gamma\) to mean that there exist increasing functions \(f\) and \(g\) such that, given \(\gamma\), any choice of \(\beta\le f(\gamma)\) and \( \alpha\le f(\beta) \) ensures the validity of all calculations in the proof.
    Moreover, when we have variable of form $1/C$ in the hierarchy, $C$ is always assumed to be a positive integer.
    
    
\begin{lemma}[\cite{han2020complexity}]\label{lem:good_parition_F}
Let $k,r,n_0\ge 2$ be integers and $\gamma,\beta >0$ where $ 1/n_0\ll \beta\ll\gamma\ll 1/r,1/k$. 
Let $F$ be an unbalanced $r$-vertex $k$-chromatic graph and $h:=r^{k-1}$. 
For each $n$-vertex graph $G$ with \(n\ge n_0\) and $\delta(G)\ge (1-1/{\chi}_{cr}(F)+\gamma)n$, in time $O(n^{2^{h-1}r+1})$ we can find an $(F,\beta,2^{h-1},1/r+\gamma/2)$-good partition $\mathcal{P}$ of $ V(G)$.
\end{lemma}

Now we are ready to state our result for $F$-factors. 
For a graph \(G\) with at least two vertices, define \(d_1(G):=e_G/(v_G-1)\) and \(m_1(G):=\max_{G'\subseteq G:v_{G'}>1}d_1(G')\). 
A graph \(F\) is \emph{strictly \(1\)-balanced} if \(d_1(F')<d_1(F)\) for every \(F' \subsetneq F\). 
The same definition will be used for hypergraphs as well.

\begin{theorem}
    [Robustness of $F$-factors]
    \label{thm:main_thm_factor}
		Let \(k,r\in \mathbb{N}\) and let \(F\) be an unbalanced $r$-vertex \(k\)-chromatic graph. Suppose that \[1/n_0\ll \beta,\mu\ll \gamma\ll 1/r,1/k.\]
        Then there exists $ C=C_{\ref{thm:main_thm_factor}}(k,r,\gamma)>0 $ such that the following holds. 
        Let $G$ be a graph on $n\ge n_0$ vertices satisfying $\delta(G)\ge (1-1/{\chi}_{cr}(F)+\gamma )n$ and \(r\mid n\).
        Let \(\mathcal{P}\) be the partition obtained from Lemma~\ref{lem:good_parition_F}. 
        If there exists an $F$-packing $ M_0 $ of size at most $(2r-1)^r$ such that $ \vec{i}_{\mathcal{P}}(V(G)\setminus V(M_0))\in L_{\mathcal{P},F}^{\mu}(G) $, and let $ G':=G-V(M_0) $, 
        then a.a.s.~$G_p'$ contains an $F$-factor provided that \(p\ge Cn^{-1/m_1(F)}\log n\).
        In particular, if \(F\) is strictly \(1\)-balanced, then a.a.s.~$G_p'$ contains an $F$-factor provided \(p\ge Cn^{-1/m_1(F)}{\log}^{1/e_F} n\).
        If no such \(F\)-packing \(M_0\) exists, then \(G\) contains no \(F\)-factor.
\end{theorem}

A special case of Theorem \ref{thm:main_thm_factor} is when $M_0=\emptyset$, and then it says that the $p$-random sparsification of $G$ a.a.s.~contains an $F$-factor (see Theorem~\ref{thm:main_technical_theorem} for a stronger result).
For the general case, obviously the existence of such $M_0$ can be tested in time $O(n^{r(2r-1)^r})$.
The difficulty for upgrading Theorem \ref{thm:main_thm_factor} to a result similar to Theorem \ref{thm:algorithm_PM} is that checking the existence of such $M_0$ may destroy the randomness of $G_p$ and we shall expand on this in the remark section. 

We remark that although Theorem~\ref{thm:main_thm_factor} is stated for unbalanced graphs \(F\), 
our framework also applies to the balanced case, where the divisibility conditions are automatically satisfied and one may take \(M_0=\emptyset\). 
In this regime, however, the result follows more directly from work of Kelly, M\"{u}yesser and Pokrovskiy~\cite{kelly2024optimal}, and we do not pursue this here.

For strictly \(1\)-balanced \(F\), Johansson, Kahn, and Vu \cite{johansson2008factors} proved that the threshold for the appearance of an \(F\)-factor in the binomial random graph is \( \Theta\left(n^{-1/m_1(F)}{\log}^{1/e_F} n\right)\). 
Hence, the assumptions on $p$ in Theorem \ref{thm:main_thm_factor} are optimal up to a constant factor for strictly \(1\)-balanced $F$.

\subsection{Spread method and the main technical result}

    To prove Theorems \ref{thm:algorithm_PM} and \ref{thm:main_thm_factor}, we use the so-called \textit{spread} method, which was introduced by Talagrand and used by Frankston, Kahn, Narayanan, and Park~\cite{frankston2021thresholds} to resolve a conjecture of Talagrand. 
    A key contribution of Talagrand's conjecture is the established connection between so-called \emph{spread measure} and threshold in random graph theory. 
    As a consequence, the task of establishing a robustness result reduces to proving the existence of a probability measure with good spread over the desired substructures. 
    Roughly speaking, this means that the probability measure selects edges of \(\mathcal{H}\) such that no specific subset of edges is overly likely to appear in an edge of \(\mathcal{H}\). 
    In other words, the distribution avoids concentrating too much weight on any particular edge set. This notion is formalized as follows.

    \begin{definition}
        Let \(q\in [0,1]\). Let \(\mathcal{H}\) be a hypergraph on vertex set $V$, and let \(\mu\) be a probability distribution on \(\mathcal{H}\). We say that \(\mu\) is \(q\)-\emph{spread} if for every set \(S\subseteq  V\):
        \[\mu(S):=\mu \left(\{A\in \mathcal{H}:S\subseteq A\}\right)\le q^{\card{S}}.\]
    \end{definition}
     In our context, we are primarily concerned with such hypergraph \(\mathcal{H}\) where \(V:=E(H)\) for some host hypergraph \(H\), and \(\mathcal{H}\) denotes the collection of \(F\)-factors in \(H\).

    Now we are ready to state our main technical result, which is a general structural theorem that serves as a robust version of the main result (Theorem~\ref{thm:han_structural_theorem}) from~\cite{han2020complexity}, strengthened via the spread-based framework described earlier. 
    Let \(F\) be an \(r\)-vertex \(k\)-graph and \(D \in \mathbb{N}\). Define \(\delta(F,\ell,D)\) to be the smallest number \(\delta\) such that every \(k\)-graph \(H\) on \(n\) vertices with \(\delta_{\ell}(H) \ge (\delta + o(1))\binom{n - \ell}{k - \ell}\) contains an \(F\)-packing covering all but at most \(D\) vertices. 
    We write \(\delta(k,\ell,D)\) for \(\delta(F,\ell,D)\) when \(F\) is a single edge. 
    Let \(\mathcal P=\{V_1,\dots,V_d\}\) be a partition. 
    Define \(L_{\max}^{d}:=\{\vec v\in\mathbb Z^{d}: r\text{ divides } \sum_{i\in[d]}\vec v|_i\}.\) 
    For any lattice \(L\subseteq L_{\max}^{d}\), let \(Q(\mathcal P,L):=L_{\max}^{d}/L.\) 

\begin{theorem}[Main technical result]\label{thm:main_technical_theorem}
    Let $k,\ell\in \mathbb{N}$ with \(\ell \le k-1\) and let $F$ be an $r$-vertex $k$-graph. 
    Define $D,C'',t,n_0,d,q\in\mathbb{N}$ and $\alpha,\beta,\mu,\gamma,\eta,\varepsilon,c>0$ where \[ 1/n_0\ll\eta\ll 1/C''\ll\beta,\mu\ll\alpha,\gamma,\varepsilon,c,1/k,1/r,1/d,1/q,1/t,1/D\]
    together with \(\alpha\ll c\).
    Let $H$ be a $k$-graph on $n\ge n_0$ vertices, with \(n\in r\mathbb{N}\). Suppose that 
    \begin{enumerate}[label=$(\roman*)$]
        \item for every \(v\in V(H)\), all but at most \(\eta\binom{n-1}{\ell-1}\) sets \(S\in\binom{V(H)\setminus\{v\}}{\ell-1}\) satisfy \(\deg_H(S\cup\{v\})\ge(\delta(F,\ell,D)+\gamma)\binom{n-\ell}{k-\ell};\) \label{item:almost_perfect_degree}
        \item $\mathcal{P}=\{V_1, \dots, V_d\}$ is an $(F,\beta,t,c)$-good partition of $V(H)$;\label{item:mathcal_P}
        \item for all but at most \(\alpha n\) vertices $ v\in V(H) $, $ \card{F^{\mu}_{\mathcal{P}}(v)}\ge \varepsilon{n}^{r-1}$;\label{item:almost_robust_edge}
        \item $|Q(\mathcal{P},L_{\mathcal{P},F}^{\mu}(H))|\le q$.
    \end{enumerate}
    If $\vec{i}_{ \mathcal{P}}(V(H))\in L_{ \mathcal{P},F}^{\mu}(H)$, then there exists a $(C''/n^{1/m_1(F)})$-spread distribution on the set of $F$-factors in $H$. 
    Moreover, if \(p\ge K_{\ref{thm:FKNP}}C''\log n/n^{1/m_1(F)} \), then a.a.s.~$H_p$ contains an $F$-factor.
\end{theorem}

    At last, this randomized construction framework enables us to derive a quantitative enumeration result.
    For perfect matchings, most existing results focus on the codegree case $\ell = k-1$. 
    For a \(k\)-uniform $n$-vertex hypergraph $H$ with minimum codegree \(\delta_{k-1}(H)\ge \delta n\) for some \(\delta>1/2\), several works show that \(H\) contains at least $(\varepsilon n)^{\frac{k-1}{k}n}$ perfect matchings; see, for example, \cite{glock2021counting, kang2024perfect, pham2022toolkit}. 
    The current best result, due to Ferber, Hardiman and Mond \cite{ferber2023counting}, proves the number of perfect matchings in such \(H\) is asymptotically close to that in the complete $k$-graph up to a multiplicative factor $\left(\frac{1}{2}+o(1)\right)^{n/k}$. 
    Subsequent work \cite{kwan2026counting} extended such asymptotic counting results to minimum \(\ell\)-degree conditions for \(\ell\ge k/2\), under degree assumptions guaranteeing the existence of a perfect matching.
    Our framework further yields a counting result for perfect matchings.  
    In particular, we show that even under the significantly weaker degree condition $\delta_{\ell}(H)\ge (c^*_{k,\ell}+\gamma)\binom{n-\ell}{k-\ell}$, corresponding to the perfect fractional matching threshold \(c^*_{k,\ell}\), either \(H\) contains no perfect matching, or it contains many of them. 
    We note that our result is optimal up to the value of $\varepsilon$.

    \begin{corollary}\label{lem:counting_PM}
        For every \(k,\ell \in \mathbb{N}\) with \(1\le\ell\le k-1\) and \(\gamma > 0\), there exists $\varepsilon>0$ such that the following holds for all sufficiently large \(n\). 
        For any $k$-graph $H$ with $\delta_{\ell}(H)\ge (c^*_{k,\ell}+\gamma)\binom{n-\ell}{k-\ell}$, $H$ contains either no perfect matching or at least $(\varepsilon n)^{\frac{k-1}{k}n}$ perfect matchings.
    \end{corollary}
    
    We also obtain a corresponding result for \(F\)-factors.

    \begin{corollary}\label{lem:counting_packing}
    For every $ r\ge 2 $ and $ \gamma > 0 $, there exists $\varepsilon>0$ such that the following holds for sufficiently large \(n\). 
    For any $r$-vertex graph $F$ and any $n$-vertex graph $G$ with $\delta(G)\ge (1-1/{\chi}_{cr}(F)+\gamma )n$, $G$ contains either no $F$-factor or at least $(\varepsilon n)^{\frac{e_F}{rm_1(F)}n}$ $F$-factors.
    \end{corollary}

\subsection{Related works}
As mentioned above, the minimum degree conditions forcing $F$-factors in graphs are well understood, and the same question for hypergraphs remains a highly challenging problem, even for perfect matchings.
On the other hand, the problem for near-perfect packings is slightly easier.
Nevertheless, we collect the results we need on $\delta(F, \ell, D)$ here.
\begin{theorem}
\cite{shokoufandeh2003proof,chang2022matching}
\label{thm:deltaFD}
For every graph $F$, $\delta(F, 1, 5v_F^2)=1-1/\chi_{cr}(F)$.
For every \(k,\ell \in \mathbb{N}\) with \(1\le\ell\le k-1\), $\delta(k,\ell,2k-\ell-1)=c^*_{k,\ell}$.
\end{theorem}


Since the work of Krivelevich, Lee and Sudakov (on Hamilton cycles), robust-type results have been proved for clique factors \cite{allen2024robust, pham2022toolkit}, hypergraph perfect matchings \cite{pham2022toolkit}, bounded-degree spanning trees \cite{pham2022toolkit}, $F$-factors \cite{kelly2024optimal}, and power of Hamilton cycles \cite{joos2023robust}. We remark that all of these results focus on minimum degree conditions as sufficient conditions, while in our result the minimum degree conditions are significantly lower. 

\section{Proof ideas}

\subsection{Proof outline for Theorem \ref{thm:algorithm_PM}}

A lattice-based framework was developed in~\cite{han2017decision, keevash2013polynomial} for deciding the existence of perfect matchings in hypergraphs under minimum codegree conditions. 
A key ingredient there is a partition lemma together with the notions of robust edge-vectors and the lattices generated by them. 
Later the framework was generalized to work under the minimum $\ell$-degree conditions \cite{han2020complexity, gan2025keevash}, and to deal with perfect \(F\)-packings in graphs~\cite{han2020complexity}. 


Suppose we have established such a partition \(\mathcal P\) and the corresponding lattice \(L\).
The key idea of the lattice-based absorption method is to reduce the search problem of perfect matchings to the search of \textit{constant-sized} matchings $M$ with certain divisibility constraint with respect to the partition $\mathcal P$ and the lattice $L$, whose existence can be determined in polynomial time.
In contrast, for our problem, we need to perform a similar search in \(H_p\):
if no such matching is found in \(H_p\), then \(H_p\) contains no perfect matching.
However, if such a matching is found in \(H_p\), we cannot argue as in the dense case by just showing \(H-V(M)\) contains a perfect matching.  
The reason is that the search procedure may expose many edges of \(H_p\), and once exposed, absent edges are deterministically revealed as non-edges.
Thus, instead of working with the \(p\)-random sparsification of \(H-V(M)\), we must work with the \(p\)-random sparsification of \(H-V(M)\) with a small set of edges removed.

To control the effect of this loss, we develop Procedure \ref{alg:perfect_matching} which adopts a selective exposure strategy which is in spirit close to kernelizations.
Roughly speaking, for each vertex and each relevant edge type, we expose edges of that type only when their total number is less than $\eta n^{k-1}$; otherwise, we leave them unexposed. 
In the latter case, by standard concentration, a.a.s.~many such edges survive in $H_p$ and remain available for the subsequent construction. 
Therefore the \textit{kernel} of our search is the former type of edges.
Since these types contribute only negligibly to the local degree structure, removing all exposed edges has only a negligible effect on the relevant properties of \(H\).
We may therefore apply Theorem~\ref{thm:main_technical_theorem}, which yields a perfect matching in \(H_p\) a.a.s.

Another difficulty arises under the minimum $\ell$-degree setting, when $\ell < k-1$.
In this case, one cannot hope for a good partition as used in Theorem \ref{thm:main_technical_theorem}, due to the existence of a small amount of vertices which are reachable to few other vertices.
In \cite{gan2025keevash}, a more complicated partition lemma was established to simultaneously classify those vertices as well, resulting in a refined partition together with two classes of robust vectors, and the robust lattice used in \cite{gan2025keevash} is generated by the union of them.
Thus, unlike in the dense setting, the two classes cannot simply be carried together into the later random redistribution step. This motivates the lattice-preparation step used in our proof. 
Roughly speaking, we show that in the proof we can always choose ahead a small number of appropriate matching edges, so that we can \textit{cancel out} the second robust lattice in a later stage of the proof.
This reduces the lattice to the ones used in the codegree setting and thus allows us to apply Theorem \ref{thm:main_technical_theorem} (see Theorem \ref{thm:auxiliary_thm_perfect_matching}).

\subsection{Proof outline for Theorem \ref{thm:main_technical_theorem}}

The main tool for the proof of Theorem \ref{thm:main_technical_theorem} is the spread method.
For this we first state a celebrated result of Frankston, Kahn, Narayanan and Park \cite{frankston2021thresholds}, conjectured by Talagrand.

\begin{theorem}\cite[Theorem 1.6]{frankston2021thresholds}\label{thm:FKNP}
        There exists a constant \(K=K_{\ref{thm:FKNP}}>0\) such that the following holds. 
        Let \(\mathcal H\) be an \(\ell\)-uniform \(q\)-spread hypergraph on vertex set \(V\). 
        If \(V_p\) is obtained by including each vertex of \(V\) independently with probability \(p\ge Kq\log \ell\), then \(V_p\) a.a.s. contains an edge in \(\mathcal H\).
\end{theorem}

While Theorem~\ref{thm:FKNP} reduces the problem to constructing a probability measure with good spread, achieving this on non-complete host graphs remains technically challenging. To address this, for our work we adopt the \emph{randomized embedding algorithm} $G\hookrightarrow H$ proposed by Kelly, M\"{u}yesser, and Pokrovskiy \cite{kelly2024optimal}. The central concept in their method is the following notion of \emph{vertex spread}. 

    \begin{definition}
        Let \(X\) and \(Y\) be finite sets, and let \(\mu\) be a probability distribution over injection \(\phi : X\rightarrow Y\). For \(q\in [0,1]\), we say that \(\mu\) is \(q\)-vertex-spread if for every two sequences of distinct vertices \(x_1,\dots, x_s\in X\) and \(y_1,\dots, y_s\in Y\), \[\mu \left(\{\phi: \phi(x_i)=y_i\text{ for all } i\in [s]\}\right)\le q^s.\]
    \end{definition}

    Based on this notion, they constructed a distribution with provably good spread properties, stated in the result below. 
    A hypergraph embedding \(\phi:G\hookrightarrow H\) of a hypergraph \(G\) into a hypergraph \(H\) is an injective map \(\phi:V(G)\hookrightarrow V(H)\) that maps edges of \(G\) to edges of \(H\), so there is an embedding of \(G\) into \(H\) if and only if \(H\) contains a subgraph isomorphic to \(G\). 

\begin{proposition}\cite[Proposition 1.17]{kelly2024optimal}\label{prop:vertex_to_edge}
        For every $C,k,\Delta>0$, there exists $C'=C_{\ref{prop:vertex_to_edge}}(C,k,\Delta)>0$ such that the following holds for all sufficiently large $n$. Let $H$ and $G$ be $n$-vertex $k$-graphs. If there is a $(C/n)$-vertex-spread distribution on embeddings $G\hookrightarrow H$ and $\Delta(G)\le \Delta$, then there is a $(C'/n^{1/m_1(G)})$-spread distribution on subgraphs of $H$ which are isomorphic to $G$.
\end{proposition}

     We now briefly outline the proof of Theorem~\ref{thm:main_technical_theorem}; the full proof is presented in Section~\ref{sec:spreadness_frm_vertex_spreadness}. 
     By Proposition~\ref{prop:vertex_to_edge}, it suffices to construct a $(C/n)$-vertex-spread distribution on embeddings $G\hookrightarrow H$, where $G$ consists of \(n/r\) disjoint copies of \(F\), that is, a perfect \(F\)-packing. 
     We begin by randomly partitioning \(V(H)\) into constant-sized clusters \(U_1,\dots,U_m\) such that each cluster inherits the minimum degree condition and partition properties from \(H\). 
     This is achieved via the \emph{random redistribution} argument introduced by~\cite{kelly2024optimal}. Although the global structure of \(H\) satisfies only approximate conditions, this step ensures that each \(U_i\) has exact degree conditions and a locally good partition. 
     We then aim to apply the general structural theorem (Theorem~\ref{thm:han_structural_theorem}) from~\cite{han2020complexity} on each cluster to find a perfect $F$-packing in \(H[U_i]\). 
     However, this theorem requires, in addition to the degree and partition assumptions, certain divisibility properties, which a cluster \(U_i\) may not initially satisfy. 
     To address this, we slightly adjust each \(U_i\) by moving a few vertices from adjacent cluster to ensure the required divisibility condition. 
     Property~\ref{item:almost_robust_edge} plays a crucial role in enabling this adjustment.
     Let \(U_i'\) denote the modified cluster, and applying the structural theorem to it gives a perfect $F$-packing \(M_i\) in \(H[U_i']\). 
     Taking the union of all such packings yields a perfect \(F\)-packing in \(H\). 
     Finally, since each cluster is sufficiently small, the probability that a fixed vertex falls in any given cluster is at most \(C/n\).
     This yields a \((C/n)\)-vertex-spread distribution on embeddings \(G \hookrightarrow H\).

    The rest of the paper is organized as follows. 
    In the next section, we introduce the basic preliminaries used throughout the paper.  
    Section~\ref{sec:random_cluster} establishes the random clustering lemma, our main tool, which is then used in Section~\ref{sec:spreadness_frm_vertex_spreadness} to prove the main technical theorem (Theorem~\ref{thm:main_technical_theorem}).
    Section~\ref{sec:robustness_F} proves the robustness of \(F\)-factors (Theorem~\ref{thm:main_thm_factor}). 
    Before turning to the algorithmic part, Section~\ref{sec:lattice_preparation} records a two-lattice preparation used in the proof of the main algorithmic result. 
    We then prove Theorem~\ref{thm:algorithm_PM} in Section~\ref{sec:proof_main_result}. 
    Finally, the proofs of the two enumeration results (Corollaries~\ref{lem:counting_PM} and~\ref{lem:counting_packing}) are presented in Section~\ref{sec:enumeration}. 
    The paper concludes with some remarks and open questions.

\section{Preliminaries}

 We need the following concentration inequality due to McDiarmid \cite{mcdiarmid1989method}, whose present formulation can be found in \cite{liebenau2023asymptotic}.
 
	\begin{lemma}\cite[Lemma 6.1]{liebenau2023asymptotic}\label{lem:McDiarmid_inequality}
	    Let $c>0$ and let $f$ be a function defined on the set of subsets of some set $V$ such that $|f(V_1)-f(V_2)|\le c$ whenever $|V_1|=|V_2|=m$ and $|V_1\cap V_2|=m-1.$ Let $A$ be a uniformly random $m$-subset of $V$. Then for any $\alpha>0$ we have \[\mathbb{P}\big[|f(A)-\mathbb{E}[f(A)]|\ge \alpha c\sqrt{m}\big]\le 2\exp{(-2\alpha^2)}.\]
	\end{lemma}

    We prove the following easy consequences of McDiarmid's inequality.

    \begin{lemma}\label{lem:application_McDiarmid}
        Let $0< \gamma'<\gamma<1$ and $1/n,1/m\ll1/s,\gamma-\gamma'$. Let $V$ be an $n$-vertex set and $\mathcal S\subseteq \binom{V}{s}$ satisfying $|\mathcal S|\ge \gamma n^s$. Let $A$ be a uniformly random $m$-subset of $V$, and define \(\mathcal S' := \{S\in \mathcal S : S\subseteq A\}.\) Then \(\mathbb{P}\left[|\mathcal S'|<\gamma' m^s\right]\le 2\exp{(-m(\gamma-\gamma')^2/2)}.\)
    \end{lemma}
    \begin{proof}
             We verify the assumptions of Lemma~\ref{lem:McDiarmid_inequality}. 
             Let $f:\binom{V}{m}\rightarrow\mathbb{R}$ be defined by $f(X):=\bigl|\{S\in\mathcal S:S\subseteq X\}\bigr|$ for each $X\in \binom{V}{m}$, and set $\varepsilon=(\gamma-\gamma')/(2\gamma)<1/2$. 
             Observe that if $X_1,X_2\in \binom{V}{m}$ with $|X_1\cap X_2|=m-1$, then changing one vertex affects only those $s$-sets containing this vertex. Hence $|f(X_1)-f(X_2)|\le m^{s-1}$, so the Lipschitz condition holds with $c=m^{s-1}$. 
             For every element \(S\in \mathcal{S}\), the probability that \(S\subseteq A\) is ${\binom{n-s}{m-s}}/{\binom{n}{m}}\ge (1-\varepsilon)\frac{{m}^{s}}{n^{s}}$, where we used $1/n,1/m\ll 1/s,\gamma-\gamma'$. 
             So by linearity of expectation, we have $\mathbb{E}[f(A)]\ge \gamma(1-\varepsilon){m}^s=\frac{\gamma+\gamma'}{2}{m}^{s}$.  
             We apply Lemma \ref{lem:McDiarmid_inequality} with $c={m}^{s-1}$ and $\alpha=\sqrt{m}(\gamma-\gamma')/2$, and get $\mathbb{P}\left[|\mathcal S'|<\gamma' m^s\right]\le 2\exp{\left(-m(\gamma-\gamma')^2/2\right)}$, as desired.
        \end{proof}

        Similarly, we obtain the following variant.
        \begin{lemma}\label{lem:application_McDiarmid_binom}
        Let $0< \gamma'<\gamma<1$ and $1/n,1/m\ll 1/s,\gamma-\gamma'$. Let $V$ be an $n$-vertex set and $\mathcal S\subseteq \binom{V}{s}$ satisfying $|\mathcal S|\ge \gamma \binom{n}{s}$. Let $A$ be a uniformly random $m$-subset of $V$, and define \(\mathcal S' := \{S\in \mathcal S : S\subseteq A\}.\) Then \(\mathbb{P}\left[|\mathcal S'|<\gamma' \binom{m}{s}\right]\le 2\exp{(-m(\gamma-\gamma')^2/4s^2)}.\)
        \end{lemma}

        The following lemmas are immediate consequences of Lemma~\ref{lem:application_McDiarmid}.

\begin{lemma}\label{lem:reachable_set}
    Let $0<\beta'<\beta<1$ and $1/n,1/m\ll 1/r,\beta-\beta'$. Let $F$ be an $r$-vertex $k$-graph, and let $H$ be an $n$-vertex $k$-graph with vertex set $V$. Suppose that the vertex set $U\subseteq V$ is $(F,\beta,t)$-closed in $H$. Let $A$ be a uniformly random $m$-subset of $V$. Then for all vertices $u,v\in U$, \[\mathbb{P}\big[u\text{ and }v\text{ are }(F,\beta',t)\text{-reachable in }H[A\cup \{u,v\}]\big]\ge 1-2\exp{\left(-m(\beta-\beta')^2/2\right).}\]
\end{lemma}

\begin{lemma}\label{lem:all_robust_edge}
    Let $0< \gamma'<\gamma<1$ and $1/n,1/m\ll 1/r,\gamma-\gamma'$. Let $F$ be an $r$-vertex $k$-graph, and let $H$ be an $n$-vertex $k$-graph with vertex set $V$ and a partition $\mathcal{P}$ of $V$. Suppose that the number of copies of $F$ with index vector $\vec{v}$ is at least $\gamma n^r $. 
    Let $A$ be a uniformly random $m$-subset of $V$. 
    Then with probability at most $2\exp{(-m(\gamma-\gamma')^2/2)}$, $H[A]$ contains fewer than $\gamma' m^r $ copies of $F$ with index vector $\vec{v}$.
\end{lemma}

For a vertex \(v\in V(H)\), a subset $Y\subseteq V(H)$, and a partition $\mathcal{P}$ of $V(H)$, we define $F^{\mu}_{\mathcal{P}}(v, Y):=\{R\in F^{\mu}_{\mathcal{P}}(v):R\setminus \{v\}\subseteq Y\}$ to be the collection of \(r\)-sets that span at least one \(\mu\)-robust copy of \(F\), contain $ v $, and have all other vertices lying entirely within \(Y\). 

\begin{lemma}\label{lem:robust_edge}
    Let $0< \gamma'<\gamma<1$ and $1/n,1/m\ll 1/r,\gamma-\gamma'$. Let $F$ be an $r$-vertex $k$-graph, and let $H$ be an $n$-vertex $k$-graph with vertex set $V$ and a partition $\mathcal{P}$ of $V$ satisfying $ \card{F^{\mu}_{\mathcal{P}}(v)}\ge \gamma n^{r-1} $ for each $v\in V$.
    Let $A$ be a uniformly random $m$-subset of $V$. 
    Then for every $v\in V$, we have \[\mathbb{P}\left[|F^{\mu}_{\mathcal{P}}(v,A)|<\gamma'{m}^{r-1 }\right]\le 2\exp{\left(-m(\gamma-\gamma')^2/2\right)}.\]
\end{lemma}

As a consequence of Lemma~\ref{lem:application_McDiarmid_binom}, we obtain the following.
\begin{lemma}\label{lem:l_degree}
    Let $0<\gamma'<\gamma<1$ and $1/n,1/m\ll 1/k,\gamma-\gamma'$. Let $H$ be an $n$-vertex $k$-graph with vertex set $V$. Suppose that \(\deg_H(D)\ge \gamma\binom{n-\ell}{k-\ell}\) for some \(\ell\)-set \(D\).
    Let $A$ be a uniformly random $m$-subset of $V$. Then \[\mathbb{P}\big[\deg_H(D,A)<\gamma'\binom{m-\ell}{k-\ell}\big]\le 2\exp{\left(-m(\gamma-\gamma')^2/4k^2\right),}\]
    where \(\deg_H(D,A)\) denotes the number of edges containing \(D\) whose remaining vertices lie in \(A\).
\end{lemma}

\begin{lemma}\cite[Lemma 3.4]{ferber2022dirac}\label{lem:almost_perfect_codegree}
    Let $0<\gamma'<\gamma<1$. Consider an $n$-vertex $k$-graph G where all but $\alpha\binom{n}{d}$ of the $d$-sets have degree at least $(\delta+\gamma)\binom{n-d}{k-d}$. Let $A$ be a uniformly random subset of $m\ge 2d$ vertices of $G$. Then with probability at least $1-\binom{m}{d}\left(\alpha+e^{-m(\gamma-\gamma')^2/4k^2}\right)$, the random induced subgraph $G[A]$ has minimum $d$-degree at least $(\delta+\gamma')\binom{m-d}{k-d}$.
\end{lemma}

We will also use a result due to McDiarmid, as presented in Chapter 16.2 of the textbook by Molloy and Reed~\cite{molloy2002graph}. In this context, a \emph{choice} refers to the position that a particular element is mapped to in a permutation.

\begin{lemma}[McDiarmid's inequality for random permutations]\label{lem:McDiarmid's inequality for random permutations}
	    Let X be a non-negative random variable determined by a random permutation $\pi$ of $[n]$ such that the following holds for some $c,r>0$:
        \begin{enumerate}
            \item Interchanging two elements of $\pi$ can affect the value of $X$ by at most $c$.
            \item For any $s$, if $X\ge s$ then there is a set of at most $rs$ choices whose outcomes certify that $X\ge s.$
        \end{enumerate}
        Then, for any $0\le t\le \mathbb{E}[X]$,
        \[
        \mathbb{P}\big[|X-\mathbb{E}[X]|\ge t+60c\sqrt{r\mathbb{E}[X]}\big]\le 4\exp{(-t^2/(8c^2r\mathbb{E}[X]))}.
        \]
\end{lemma}

Finally, we need the following lemma, recorded in
\cite{kelly2024optimal} and based on the random sampling approach from \cite{pham2022toolkit}.
        
\begin{lemma}\cite[Lemma 3.2]{kelly2024optimal}\label{lem:spread_PM}
    There exists an absolute constant $C_{\ref{lem:spread_PM}}$ with the following property. If $G$ is a balanced bipartite graph on $2n$ vertices with $\delta(G)\ge 3n/4$, then there exists a $(C_{\ref{lem:spread_PM}}/n)$-spread distribution on perfect matchings of $G$.
\end{lemma}

\section{Random clustering lemma}\label{sec:random_cluster}

Our main tool for the proof of Theorem~\ref{thm:main_technical_theorem} is Lemma~\ref{lem:random_cluster_degree}, which is proved by the random clustering method of Kelly, M\"{u}yesser and Pokrovskiy~\cite{kelly2024optimal}, and we use it to construct a probability distribution over the desired subgraphs with good spread. For $Y\subseteq V(H)$ and a partition $\mathcal{P}$ of $V(H)$, we define the projection of $\mathcal{P}$ on $Y$ as  $\mathcal{P}^Y:=\{W\cap Y:W\in \mathcal{P}\}$.
    
Now we are ready to state the main result of this section, which says that we can partition the vertex set of a dense $k$-graph $H$ randomly, so that all parts inherit the properties of $H$, and in particular, the probability that a vertex is in a specific part is at most $O(1/n)$.
Note that to make it applicable in the proof of Theorem~\ref{thm:main_technical_theorem} (assumption \ref{item:almost_perfect_degree} therein), we do not require the degree condition to hold for all \((\ell-1)\)-sets; instead, for each \(v\in V(H)\), it may fail for up to an \(\eta\)-fraction of such sets. 

\begin{lemma}[Random clustering lemma for \(\ell\)-degree]\label{lem:random_cluster_degree}
Let $k,\ell\in \mathbb{N}$ and let $F$ be an $r$-vertex $k$-graph.
Let $d,q,t,C,n\in \mathbb{N}$ and $\alpha,\beta,\mu,\gamma,\eta,\varepsilon,c>0$ satisfy 
\[1/n\ll\eta\ll1/C'\ll 1/C\ll \alpha,\beta,\mu,\gamma,\varepsilon,c,1/k,1/r,1/d,1/q \le 1,\]
together with \(\alpha\ll c\), $r\mid n$,  and $r\mid C$. Let $H$ be an $n$-vertex $k$-graph. Suppose that:
\begin{enumerate}[label=(P\arabic*)]
\item for every \(v\in V(H)\),  all but at most \(\eta\binom{n-1}{\ell-1}\) sets \(S\in\binom{V(H)\setminus\{v\}}{\ell-1}\) satisfy \(\deg_H(S\cup\{v\})\ge(\delta+\gamma)\binom{n-\ell}{k-\ell};\)\label{condition:almost_perfect_degree}
\item $\mathcal{P}=\{V_1, \dots, V_d\}$ is an $(F,\beta,t,c)$-good partition of $V(H)$;\label{condition:good_partition}
\item for all but at most \(\alpha n\) vertices $ v\in V(H) $, $\card{F^{\mu}_{\mathcal{P}}(v)}\ge \varepsilon{n}^{r-1}$.\label{condition:robust_edges}
\end{enumerate}
Then there exists a random partition $\mathcal{U}=\{U_1,U_2,\dots,U_m\}$ of $V(H)$ with the following properties:
\begin{enumerate}[label=(R\arabic*)]
\item $\card{U_1}$ equals $(C-1)C$ plus the remainder when $n$ is divided by $(C-1)C$, and $|U_i|=C$ for every $2\le i\le m$;\label{item:(R1)}

\item for every $i\in[m]$, \(\delta_{\ell}(H[U_i])\ge (\delta+\gamma/2)\binom{\card{U_i}-\ell}{k-\ell}\), and \(\mathcal{P}^{U_i}\) is an \((F,\beta/2,t,c/2)\)-good partition of \(H[U_i]\);\label{item:(R-degree)}

\item for every $2\le i\le m$, there exists a set $T_i\subseteq U_i$ with $|T_i\cap V_j|=rq$ for each $j\in[d]$ such that every vertex $v\in T_{i}$ satisfies $\card{F^{\mu}_{\mathcal{P}}(v, U_{i-1})}\ge \varepsilon{\card{U_{i-1}}}^{r-1}/2 $;\label{item:(R-robust)}

\item for every $i\in [m]$, 
$I_{\mathcal{P},F}^{\mu}(H)\subseteq I_{ \mathcal{P}^{U_i},F}^{\mu/2}(H[U_i])$;\label{item:(R4)}
            
\item for every integer \(s\in [n]\), every set of distinct vertices $y_1,\dots,y_s\in V(H)$, and every function $f:[s]\rightarrow[m]$, \label{item:(R5)}
\[\mathbb{P}\big[y_i\in U_{f(i)}\text{ for each }i\in[s]\big]\le \Big({\frac{C'}{n}}\Big)^s. 
\]
\end{enumerate}
\end{lemma}

Now let us explain the properties of the random partition $U_i$ found by Lemma~\ref{lem:random_cluster_degree}.
First,~\ref{item:(R5)} reflects the probability of the location of the vertices, which is exactly what we want for getting the vertex-spread.
The other properties together are designed so that we could find an $F$-factor in each cluster $U_i$: the size of the clusters must be a multiple of $r$, every part must inherit the minimum degree of $H$, and the reachability property and the robust-lattice structure should also be inherited.
Moreover, in~\ref{item:(R-robust)} we include small sets $T_i\subseteq U_i$, which will be used to correct the divisibility in the proof of Theorem~\ref{thm:main_technical_theorem} -- when $\vec{i}_{\mathcal P^{U_i}}(U_i)$ is not in $L_{\mathcal P^{U_i}}^{\mu/2}(H[U_i])$, 
we will add a subset of $T_{i+1}$ to $U_i$ so that the new index vector does belong to the lattice, which allows us to build an $F$-factor in the $i$-th cluster.
Then this correction can be done in a recursive manner.

Before presenting the full proof, we outline the key ideas behind it. The argument is based on the random redistribution technique introduced in~\cite{kelly2024optimal}, and we modify their framework to suit our setting by adjusting both the technical assumptions and the properties we aim to achieve. 
We begin with randomly partitioning \(V(H)\) into clusters \(U_1',\dots,U_m'\). 
This initial partition is, in a sense, already ``close" to the one required in Lemma~\ref{lem:random_cluster_degree}: the fifth property is satisfied, and the remaining properties are violated by only a small number of clusters. 
We refer to a cluster \(U_i'\) as \emph{bad} if it violates any of the required properties stated in the lemma. 
Following the idea of~\cite{kelly2024optimal}, we then randomly redistribute the vertices from the bad clusters into the good ones in a way that preserves the key structural properties of the latter. 
    
For the remainder of this section, fix 
\[1/n\ll\eta\ll1/C'\ll 1/C\ll \alpha,\beta,\mu,\gamma,\varepsilon,c,1/k,1/r,1/d,1/q \le 1,\]
together with \(\alpha\ll c\), $r\mid n$,  and $r\mid C$, as in Lemma \ref{lem:random_cluster_degree}. 
Let $C_1$ be $(C-1)C$ plus the remainder when $n$ is divided by $(C-1)C$. Let $W_1:=[C_1]$, and for $i\ge2$, let $W_i:=[C_1+(C-1)(i-1)]\setminus{\bigcup_{j=1}^{i-1}W_j}=[C_1+(C-1)(i-2)+1,C_1+(C-1)(i-1)].$ 
Note that $\{W_1,\dots,W_m\}$ is a partition of $[n]$ for $m:=(n-C_1)/(C-1)+1$. 
Let $H$ be a $k$-graph with vertex set $V:=\{v_1,\dots,v_n\}$ 
and let $\mathcal{P}=\{V_1,\dots,V_d\}$ be a partition of $V$ as described in Lemma~\ref{lem:random_cluster_degree}. 
Let $\pi:[n]\rightarrow[n]$ be a uniformly random permutation of $[n]$, and let $U_i':=\{v_j:\pi(j)\in W_i\}.$ 
    
To prepare for the proof of Lemma~\ref{lem:random_cluster_degree}, we begin by establishing a sequence of auxiliary lemmas, all under the assumptions of Lemma~\ref{lem:random_cluster_degree}. 
Each lemma verifies a specific aspect of the desired structure, and together they provide the framework needed for the random redistribution argument to succeed. 
While these lemmas address the individual components, the complete proof will be presented at the end of this section.

The following lemma bounds the probability that Property~\ref{item:(R-degree)} holds. 
Let $ C^*:=|U_i'|$, that is, $C^*=C_1$ if $i=1$ and $C^*=C-1$ if $i>1$.
So $C^*\ge C-1\ge C/2$. 
Note that Lemmas~\ref{lem:prob_perfect_condition}--\ref{lem:U_1} are under the assumptions of Lemma~\ref{lem:random_cluster_degree}.

\begin{lemma}\label{lem:prob_perfect_condition}
    Let \(\theta=\min\{\beta,c,\gamma/k^2\}\). For every vertex \(v\in V\) and every \(i\in [m]\), 
    \begin{equation}\label{eq:l_degree}
        \mathbb{P}\left[\delta_{\ell}(H[U_i'\cup \{v\}])<(\delta+2\gamma/3)\binom{\card{U_i'\cup \{v\}}-\ell}{k-\ell}\right]\le \exp{(-\theta^2C/100)},
    \end{equation} 
    \begin{equation}\label{eq:reachable_set}
        \mathbb{P}\left[\exists j\in[d],V_j\cap (U_i'\cup \{v\})\text{ is not }(F,2\beta/3,t)\text{-closed in }H[U_i'\cup \{v\}]\right]\le \exp{(-\theta^2C/100)},
    \end{equation} and 
    \begin{equation}\label{eq:size}
        \mathbb{P}\left[\exists j\in[d],|V_j\cap (U_i'\cup \{v\})|<2c\card{U_i'\cup \{v\}}/3\right]\le \exp{(-\theta^2C/100)}.
    \end{equation}
\end{lemma}

\begin{proof}

     Fix \(v\in V\) and \(i\in [m]\). Call an \((\ell-1)\)-set \emph{bad} for \(v\) if \(\deg_H(S\cup \{v\})< (\delta+\gamma)\binom{n-\ell}{k-\ell}\). Note that the number of such bad \((\ell-1)\) sets in \(H\) is at most \(\eta\binom{n-1}{\ell-1}\), that  $V_j$ is $(F,\beta,t)$-closed in $H$, and that $|V_j|\ge cn$ for all $j\in[d]$.
    
    \medskip
    \noindent
    \textbf{Proof of \eqref{eq:reachable_set}.} 
    Fix \(j\in [d]\), let $A_j$ be the event that there exists a pair of vertices in $V_{j}\cap (U_i'\cup\{v\})$ that are not \((F,2\beta/3,t)\)-reachable in \(H[U_i']\).
    By Lemma~\ref{lem:reachable_set} with \(V_j\) playing the role of \(U\) and \(U_i'\) playing the role of \(A\), we have 
    \[
    \mathbb{P}\left[A_j\right]\le \binom{C^*+1}{2}\cdot 2\exp{(-\beta^2C^*/32)}\le \exp{(-\theta^2C/50)}.
    \] 
    Since \(C\) is sufficiently large, every pair of vertices in \(V_j\cap U_i'\) that are \((F,3\beta/4,t)\)-reachable in \(H[U_i']\) remain \((F,2\beta/3,t)\)-reachable in \(H[U_i'\cup\{v\}]\). 
    Therefore, 
    \[\mathbb{P}\left[V_j\cap (U_i'\cup\{v\}) \text{ is not }(F,2\beta/3,t)\text{-closed in }H[U_i'\cup\{v\}]\right]\le \exp{(-\theta^2C/50)}.\] 
    Taking a union bound over \(j\in[d]\), we obtain \eqref{eq:reachable_set}.

    \medskip
        \noindent
\textbf{Proof of \eqref{eq:size}.} 
    For every \(j\in [d]\), applying Lemma~\ref{lem:application_McDiarmid} with $V_j$ as $\mathcal{S}$ and $U_i'$ as $A$ gives
    \[
    \mathbb{P}\Big[|V_j\cap U_i'|< 3c|U_i'|/4\Big] \le 2 \exp(-c^2 C^*/36) \le \exp(-\theta^2 C/50).
    \] 
    Moreover, we have \(\mathbb{P}\Big[|V_j\cap (U_i'\cup\{v\})|< 2c|U_i'\cup\{v\}|/3\Big]\le \mathbb{P}\Big[|V_j\cap U_i'|< 3c|U_i'|/4\Big]\) since \(C\) is sufficiently large. 
    A union bound over $j\in[d]$ gives \eqref{eq:size}.
    
    \medskip
        \noindent
\textbf{Proof of \eqref{eq:l_degree}.} 
    We consider $\ell$-sets in $H[U_i'\cup\{v\}]$ in two types: those containing $v$ and those not containing $v$. 

    \textbf{Case 1: $\ell$-sets containing \(v\).}

    On the one hand, let \(X_{v,i}\) be the number of bad \((\ell-1)\)-sets contained in \(U_i'\).
    Using linearity of expectation together with \ref{condition:almost_perfect_degree}, we have \(\mathbb{E}\left[X_{v,i}\right]\le \eta \binom{n-1}{\ell -1}{\binom{n-(\ell-1)}{\card{U_i'}-(\ell-1)}}/{\binom{n}{\card{U_i'}}}\le \eta {\card{U_i'}}^{\ell-1}\). 
    Thus, by Markov's inequality, \(\mathbb{P}[X_{v,i}\ge 1] \le \eta {\card{U_i'}}^{\ell-1}\). 
    On the other hand, fix an $(\ell-1)$-set $S\subseteq U_i'$ that is not bad for $v$. 
    Lemma~\ref{lem:l_degree} then gives \(\mathbb{P}\left[\deg_H(S\cup \{v\}, U_i')<(\delta+3\gamma/4)\binom{|U_i'|-\ell}{k-\ell}\right]\le 2\exp(-\gamma^2 C^*/32k^2).\)
    Observe that the probability that there exists an $(\ell-1)$-set $S\subseteq U_i'$ with \(\deg_H(S\cup \{v\}, U_i')<(\delta+3\gamma/4)\binom{|U_i'|-\ell}{k-\ell}\) is bounded by the probability that either:
\begin{itemize}
\item there exists a bad $(\ell-1)$-set in $U_i'$ or
\item there exists a non-bad $(\ell-1)$-set that fails the required degree condition in $U_i'$.
\end{itemize}
    Hence, this probability is at most \(\eta {\card{U_i'}}^{\ell-1}+\binom{|U_i'|}{\ell-1}\cdot 2\exp(-\gamma^2 C^*/32k^2)\le \exp{(-\theta^2C/50)}\), as $\eta\ll 1/C, \theta$.

    \textbf{Case 2: $\ell$-sets not containing \(v\).}
    
    These are precisely the $\ell$-sets contained in $U_i'$. Note that all but $\eta\binom{n}{\ell}$ of the $\ell$-sets in \(H\) have degree at least $(\delta+\gamma)\binom{n-\ell}{k-\ell}$. Applying Lemma~\ref{lem:almost_perfect_codegree}, we have \[\mathbb{P}\left[\delta_{\ell}(H[U_i'])<(\delta+3\gamma/4)\binom{\card{U_i'}-\ell}{k-\ell}\right]\le \binom{C^*}{\ell}\left(\eta+\exp{(-\gamma^2C^*/64k^2)}\right)\le \exp{(-\theta^2 C/80)}.\]

    Finally, applying a union bound over the two types of $\ell$-sets and using that \(C\) is sufficiently large, we obtain \[\mathbb{P}\left[\delta_{\ell}(H[U_i'\cup \{v\}])<(\delta+2\gamma/3)\binom{\card{U_i'\cup \{v\}}-\ell}{k-\ell}\right]\le \exp{\left(-\frac{\theta^2C}{50}\right)}+\exp{\left(-\frac{\theta^2C}{80}\right)}\le \exp{\left(-\frac{\theta^2C}{100}\right)}.\]

\end{proof}

Lemma~\ref{lem:prob_perfect_condition} immediately implies the following lemma, which guarantees that property~\ref{item:(R-degree)} holds for the first cluster.

\begin{lemma}\label{lem:U_1_condition}
    Let \(E_1\) be the event that \(\delta_{\ell}(H[U_1'])\ge(\delta+2\gamma/3)\binom{\card{U_1'}-\ell}{k-\ell}\) and \(\mathcal{P}^{U_1'}\) is an \((F,2\beta/3,t,c/2)\)-good partition of \(H[U_1']\). Then $\mathbb{P}\left[E_1\right]\ge 99/100$.
\end{lemma}

To ensure property~\ref{item:(R-degree)} holds for the remaining clusters, we next establish two lemmas that address its key components.  
Lemma~\ref{lem:d(v)} ensures that, at the random redistribution stage, each vertex has many good clusters into which it can be added.

\begin{lemma}\label{lem:d(v)}
    Let $E_2$ be the event that for every vertex $v\in V$, there are at least $(1-1/C^2)m$ indices \(i\in \{2,\dots,m\}\) such that \(\delta_{\ell}(H[U_i'\cup \{v\}])\ge(\delta+2\gamma/3)\binom{\card{U_i'\cup \{v\}}-\ell}{k-\ell}\) and  \(\mathcal{P}^{U_i'\cup \{v\}}\) is an \((F,2\beta/3,t,2c/3)\)-good partition of \(H[U_i'\cup \{v\}]\).
    Then, $\mathbb{P}\left[E_2\right]\ge 99/100$.
\end{lemma}
    
\begin{proof}
For every \(v\in V\), let $X_v$ denote the number of random sets $U_i'$ such that \(\delta_{\ell}(H[U_i'\cup \{v\}])\ge (\delta+2\gamma/3)\binom{\card{U_i'\cup \{v\}}-\ell}{k-\ell}\), \(|V_j\cap (U_i'\cup \{v\})|\ge 2c\card{U_i'\cup \{v\}}/3\) and $V_j\cap (U_i'\cup\{v\})$ is $(F,2\beta/3,t)$-closed in $H[U_i'\cup \{v\}]$ for all $j\in[d]$. 
By linearity of expectation and Lemma~\ref{lem:prob_perfect_condition}, we have \[\mathbb{E}[X_v]\ge \left(1-3\exp{(-\theta^2C/100)}\right)m.\] 
Note that interchanging two elements of $\pi$ can affect the value of $X_v$ by at most 2 and if $X_v\ge x$, this can be certified by at most $2C^2x$ choices of the random permutation. 
Therefore, by Lemma \ref{lem:McDiarmid's inequality for random permutations} applied with $c=2$, $r=2C^2$, and $t=m/C^3$, we have
\[\mathbb{P}\big[X_v<(1-1/C^2)m\big]\le 4\exp{\Big(-\frac{{(m/C^3)}^2}{64C^2m}\Big)}\le \exp{(-n/C^{10})}.\]
By the union bound, we have that $X_v\ge (1-1/C^2)m$ for every vertex $v\in V$ with probability at least 99/100.
\end{proof}

Lemma~\ref{lem:d(i)} guarantees that most clusters admit many vertices that can be added while preserving the desired properties.

\begin{lemma}\label{lem:d(i)}
    Let $E_3$ be the event that the following hold for all but at most 
    \(\exp{\left(-\theta^2 C/500\right)}m\) indices \(i\in \{2,\dots,m\}\), there exists at least \((1-\exp{\left(-\theta^2 C/200\right)})n\) vertices \(v\in V\) such that $\delta_{\ell}(H[U_i'\cup \{v\}])\ge (\delta+2\gamma/3)\binom{|U_i'\cup\{v\}|-\ell}{k-\ell}$ and \(\mathcal{P}^{U_i'\cup \{v\}}\) is an \((F,2\beta/3,t,2c/3)\)-good partition of \(H[U_i'\cup \{v\}]\). 
    Then, $\mathbb{P}\left[E_3\right]\ge 99/100$.
\end{lemma}

\begin{proof}
Call a vertex $v$ \emph{bad} for $U_i'$ if $\delta_{\ell}(H[U_i'\cup \{v\}])< (\delta+2\gamma/3)\binom{|U_i'\cup\{v\}|-\ell}{k-\ell}$ holds or \(\mathcal{P}^{U_i'\cup \{v\}}\) is not an \((F,2\beta/3,t,2c/3)\)-good partition of \(H[U_i'\cup \{v\}]\). 
For any $U_i'$, by Lemma \ref{lem:prob_perfect_condition} and linearity of expectation, we have $\mathbb{E}\big[|\{v\in V: v\text{ is bad for }U_i'\}|\big]\le 3\exp{\left(-\theta^2 C/100\right)}n$. 
Thus, by Markov's inequality, $\mathbb{P}\big[|\{v\in V: v\text{ is bad for }U_i'\}|\ge \exp{\left(-\theta^2 C/200\right)}n\big]\le 3\exp{\left(-\theta^2 C/200\right)}$. 
Let $X$ be the number of $i\in\{2,\dots,m\}$ such that $|\{v\in V: v\text{ is bad for }U_i'\}|\ge \exp{\left(-\theta^2 C/200\right)}n$. 
By linearity of expectation, $\mathbb{E}[X]\le 3\exp{\left(-\theta^2 C/200\right)}m$. 
Applying Markov's inequality again, we obtain $\mathbb{P}[X\ge\exp{\left(-\theta^2 C/500\right)}m]\le 3\exp{\left(-\theta^2 C/500\right)}\le 1/100$, as desired. 
\end{proof}

The next lemma will be used in the proof of \ref{item:(R4)}, providing a key probabilistic bound. 

\begin{lemma}\label{lem:index_vector}
Let $E_4$ be the event that the following hold for at least $\left(1-\exp{(-\mu^2C/50)}\right)m$ indices $i\in\{2,\dots,m\}$ and for $i=1$:
every $\vec{v}\in I_{\mathcal{P},F}^{\mu}(H)$ satisfies $\vec{v}\in I_{\mathcal{P}^{U_i'},F}^{2\mu/3}(H[U_i'])$. 
Then, $\mathbb{P}\left[E_4\right]\ge 99/100$.
\end{lemma}
\begin{proof}
For every $\vec{v}\in I_{\mathcal{P},F}^{\mu}(H)$ and $i\in[m]$, we say that $\vec{v}$ is \emph{bad} for $U_i'$ if $H[U_i']$ contains fewer than $2\mu {|U_i'|}^r/3$ copies of $F$ with index vector $\vec{v}$. 
Since \( \mathcal{P}^{U_i'} \) is the restriction of \( \mathcal{P} \) to \( U_i' \), we have \(\vec{i}_{\mathcal{P}^{U_i'}}(S)=\vec{i}_{\mathcal{P}}(S)\) for any subset \(S\subseteq U_i'\). 
It follows that if $\vec{v}$ is not bad for $U_i'$, then $\vec{v}\in I_{\mathcal{P}^{U_i'},F}^{2\mu/3}(H[U_i'])$. 

Next, we estimate the probability that a fixed $\vec{v}\in I_{\mathcal{P},F}^{\mu}(H)$ is bad for some $U_i'$; by Lemma \ref{lem:all_robust_edge} (with $U_i'$ playing the role of $A$), this probability is at most $2\exp{(-\mu^2C/20)}$. 
Note that $\cardi{I_{\mathcal{P},F}^{\mu}(H)}\le \binom{r+d-1}{r}$.
Then, by the union bound, the probability that there exists some \(\vec{v}\in I_{\mathcal{P},F}^{\mu}(H)\) bad for \(U_1'\) is at most \[2\exp{(-\mu^2C/20)}\binom{r+d-1}{r}\le \exp{(-\mu^2C/50)}.\] 
On the other hand, let $X$ denote the number of indices \(i\in \{2,\dots,m\}\) for which there exists some $\vec{v}\in I_{\mathcal{P},F}^{\mu}(H)$ that is bad for $U_i'$. 
By linearity of expectation, $\mathbb{E}[X]\le 2\exp{(-\mu^2C/20)}\binom{r+d-1}{r}m.$ 
By Markov's inequality, \[\mathbb{P}\big[X\ge\exp{(-\mu^2C/50)}m\big]\le 2\exp{(-\mu^2C/50)}\binom{r+d-1}{r}\le \exp{(-\mu^2C/100)}.\] 
Combining this with the bound for \(U_1'\), and applying the union bound, we obtain the desired conclusion.
\end{proof}

The following two lemmas are crucial for establishing conclusion \ref{item:(R-robust)} of the theorem. 
We call a vertex \(v\in V\) \emph{good} if $\card{F^{\mu}_{\mathcal{P}}(v)}\ge \varepsilon{n}^{r-1}$; otherwise, \(v\) is \emph{bad}. For every $i\in [m]$, let \(B_i\) be the number of bad vertices in \(U_i'\). 

\begin{lemma}\label{lem:number_bad_vertices}
Let $E_5$ be the event that the following hold: 
\begin{itemize}
\item For all but at most  $\exp{(-\alpha^2C/20)}m$ indices $i\in\{2,\dots,m\}$, $B_i < 2\alpha\card{U_i'}$;
\item The same holds for \(U_1'\), i.e., 
$B_1 < 2\alpha\card{U_1'}$.
\end{itemize} 
Then, $\mathbb{P}\left[E_5\right]\ge 99/100$.
\end{lemma}
\begin{proof}
For every $i\in [m]$, by Lemma~\ref{lem:application_McDiarmid} and \ref{condition:robust_edges}, we have \[\mathbb{P}\left[B_i\ge 2\alpha|U_i'|\right]\le 2\exp{(-\alpha^2C^*/2)}\le \exp{(-\alpha^2C/10)}.\]
Let $X$ denote the number of indices \(i\in \{2,\dots,m\}\) such that $B_i\geq 2\alpha|U_i'|$. 
By linearity of expectation, $\mathbb{E}[X]\le \exp{(-\alpha^2C/10)}m.$ 
Applying Markov's inequality, we have $\mathbb{P}\big[X\ge\exp{(-\alpha^2C/20)}m\big]\le \exp{(-\alpha^2C/20)}.$ 
Together with the probability bound for \(U_1'\), we obtain the desired conclusion.
\end{proof}

\begin{lemma}\label{lem:d_D(v)}
    Let $E_6$ be the event that for every good vertex $v\in V$, there are at least \((1-1/C^2)m\) indices \(i\in \{2,\dots,m\}\) such that $\card{F^{\mu}_{\mathcal{P}}(v, U_i')}\ge 2\varepsilon{\card{U_{i}'}}^{r-1}/3 $.
    Then, $\mathbb{P}\left[E_6\right]\ge 99/100$.
\end{lemma}

\begin{proof}
Note that for each good vertex $ v\in V(H) $, $ |F^{\mu}_{\mathcal{P}}(v)|\ge \varepsilon{n}^{r-1} $  as stated in Lemma \ref{lem:random_cluster_degree}. 
For every good vertex $v\in V$ and $i\in[m]$, by Lemma \ref{lem:robust_edge} with $U_i'$ playing the role of $A$, we have
\begin{equation}
\label{eq:FPvU}
\mathbb{P}\left[|F^{\mu}_{\mathcal{P}}(v,U_i')|<2\varepsilon{\card{U_i'}}^{r-1}/3\right]\le 2\exp{\left(-\varepsilon^2 C/20\right)}.
\end{equation}
Let $Y_v$ denote the number of $i\in [m]$ such that $|F^{\mu}_{\mathcal{P}}(v,U_i')|\ge {2\varepsilon}{\card{U_i'}}^{r-1}/3$. 
By linearity of expectation, $\mathbb{E}[Y_v]\ge \left(1-2\exp{\left(-\varepsilon^2 C/20\right)}\right)m.$ 
Note that interchanging two elements of $\pi$ can affect the value of $Y_v$ by at most 2 and if $Y_v\ge x$, this can be certified by at most $2C^2x$ choices of the random permutation. 
Therefore, by Lemma \ref{lem:McDiarmid's inequality for random permutations} applied with $c=2$, $r=2C^2$, and $t=m/C^3$, we have
\[\mathbb{P}\big[Y_v<(1-1/C^2)m\big]\le 4\exp{\Big(-\frac{{(m/C^3)}^2}{64C^2m}\Big)}\le \exp{(-n/C^{10})}.\]
By the union bound, we have that $Y_v\ge (1-1/C^2)m$ for every vertex $v\in V$ with probability at least 99/100, as desired.
\end{proof}

We shall show that $U_1'$ satisfies the properties~\ref{item:(R1)}--\ref{item:(R4)} with $U_1=U_1'$. The next lemma is important for establishing that $U_1'$ satisfies \ref{item:(R-robust)}. 
For every \(i\in [m]\), choose some \(T_i \subseteq U_i'\) of size \(rqd\) as follows. Let \(\mathcal{T}_i\) denote the collection of all \(rqd\)-subsets \(T_i\subseteq U_i'\) such that \(\card{T_i\cap V_j}=rq\) for all \(j\in [d]\), and every vertex in \(T_i\) is good. 
        \begin{itemize}
            \item If \(\mathcal{T}_i\neq\emptyset\), then we fix \(T_i\in\mathcal{T}_i\) arbitrarily.
            \item If \(\mathcal{T}_i=\emptyset\), then we choose any \(T_i\in\binom{U_i'}{rqd}\) arbitrarily.
        \end{itemize} 
        Given such a choice of \(T_1,\dots,T_m\), we consider an auxiliary digraph \(D\) with vertex set \([m]\), where \((i,j)\in E(D)\) if for every vertex \(v\in T_i\), we have \(|{F^{\mu}_{\mathcal{P}}(v, U_j')}|\ge 2\varepsilon {|{U_j'}|}^{r-1}/3.\) 
        Note that property~\ref{item:(R-robust)} in Lemma~\ref{lem:random_cluster_degree} requires that the sets \(T_i\) lie in \(\mathcal{T}_i\). Although it is possible that \(\mathcal{T}_i=\emptyset\) for some \(i\), Lemmas~\ref{lem:d(i)} and \ref{lem:number_bad_vertices} imply that the number of such ``bad'' indices is small with high probability. In the proof of Lemma~\ref{lem:random_cluster_degree}, we will reassign these bad clusters and retain only those \(U_i\) for which \(\mathcal{T}_i\neq\emptyset\), thereby ensuring that the final selection of \(T_i\) satisfies the desired property.

\begin{lemma}\label{lem:U_1}
Let $E_7$ be the event that $U_1'$ satisfies $d_D^-(1)\ge (1-\exp{\left(-\varepsilon^2 C/50\right)})(m-1)$. 
Then $\mathbb{P}\left[E_7\right]\ge 99/100$.
\end{lemma}
\begin{proof}
Recall that if \(v\in V\) is good, then~\eqref{eq:FPvU} 
says $\mathbb{P}\left[|F^{\mu}_{\mathcal{P}}(v,U_1')|<2\varepsilon{\card{U_1'}}^{r-1}/3\right]\le 2\exp{\left(-\varepsilon^2 C/20\right)}$. 
This implies that for every \(i\in \{2,\dots,m\}\), \(\mathbb{P}\left[(i,1)\in E(D)\right]\ge 1-2\exp{\left(-\varepsilon^2 C/20\right)}rqd.\)
By linearity of expectation, \(\mathbb{E}\left[d_D^-(1)\right]\ge \left(1-2\exp{\left(-\varepsilon^2 C/20\right)}rqd\right)\left(m-1\right)\). 
Then by Markov's inequality, 
\[
\mathbb{P}\left[m-1-d_D^-(1)\ge \exp{\left(-\varepsilon^2 C/50\right)\left(m-1\right)}\right]\le 2\exp{\left(-\varepsilon^2 C/50\right)}rqd,
\] implying the desired result.
\end{proof}

We now proceed to prove Lemma~\ref{lem:random_cluster_degree} by considering the distribution of \(U_1',\dots,U_m'\) conditioning on the events established in the previous lemmas, and then applying the random redistribution argument.

\begin{proof}[Proof of Lemma~\ref{lem:random_cluster_degree}]
We condition on $E_1\cap \cdots\cap E_7$, which holds with probability at least $93/100$ by Lemmas \ref{lem:U_1_condition}--\ref{lem:U_1}. 
For $i\in [m]$ such that \(T_i\in \mathcal{T}_i\), because all vertices in $T_i$ are good vertices, by Lemma~\ref{lem:d_D(v)}, we obtain that \(d_D^+(i)\ge (1-rqd/C^2)m\).

For each permutation $\pi:[n]\rightarrow [n]$, define a set of ``bad'' clusters $\mathcal{F}_{\pi}\subseteq\{U_2',\dots,U_m'\}$ which includes $U_i'$ if any of the following holds:

\begin{enumerate}[label=(A\arabic*)]
    \item $\delta_{\ell}(H[U_i'\cup \{v\}])< (\delta+2\gamma/3)\binom{|U_i'\cup \{v\}|-\ell}{k-\ell}$ or \(\mathcal{P}^{U_i'\cup \{v\}}\) is not an \((F,2\beta/3,t,2c/3)\)-good partition of \(H[U_i'\cup \{v\}]\) for at least \(\exp{\left(-\theta^2 C/200\right)}n\) vertices \(v\in V\),\label{item:bad_cluster_reachable}
    \item \(I_{\mathcal{P},F}^{\mu}(H)\nsubseteq I_{ \mathcal{P}^{U_i'},F}^{2\mu/3}(H[U_i'])\), \label{item:bad_cluster_index_verctor}
    \item \(|B_i|\ge 2\alpha |U_i'|\), and \label{item:bad_cluster_bad_vtx}
    \item $d_D^-(i)<(1-rqd/\sqrt{C})m$.\label{item:bad_cluster_T_i}
\end{enumerate}

Let $m':=(m-1)(C-1)/C$, which is an integer by the definition of \(C_1\). 
We claim that if $\pi \in E_1\cap \cdots\cap E_7$, then $|\mathcal{F}_{\pi}|\le m-1-m'=m-1-(m-1)(C-1)/C=(m-1)/C$. 
Indeed, 
there are at most $\exp{\left(-\theta^2 C/500\right)}m$ $U_i'\in \mathcal{F}_{\pi}$ of type \ref{item:bad_cluster_reachable} by Lemma \ref{lem:d(i)} assuming $\pi \in E_3$, 
at most $\exp{(-\mu^2C/50)}m$ $U_i'\in \mathcal{F}_{\pi}$ of type \ref{item:bad_cluster_index_verctor} by Lemma \ref{lem:index_vector} assuming $\pi \in E_4$,
and there are at most $\exp{\left(-\alpha^2 C/20\right)}m$ $U_i'\in \mathcal{F}_{\pi}$ of type \ref{item:bad_cluster_bad_vtx} by Lemma \ref{lem:number_bad_vertices} assuming $\pi \in E_5$.  
Next, let \(\mathcal{I}\subseteq [m]\) be the set of indices \(i\) for which \(T_i\in\mathcal{T}_i\). 
If a cluster \(U_i'\) is not of type~\ref{item:bad_cluster_bad_vtx} and satisfies \(|V_j\cap U_i'|\ge 2c|U_i'|/3\) for all \(j\in[d]\), then since at most \(2\alpha |U_i'|\) vertices in \(U_i'\) are bad and \(\alpha\ll c\), each set \(V_j\cap U_i'\) contains at least \(rq\) good vertices. 
Hence \(\mathcal{T}_i\neq\emptyset\), and by construction \(T_i\in\mathcal{T}_i\). 
Using the above bounds, this implies 
\[
|\mathcal{I}|\ge \left(1-\exp(-\theta^2 C/500)-\exp(-\alpha^2 C/20)\right)m.
\]
Note that \(\sum_{i=1}^m d_D^-(i)=\sum_{i=1}^m d_D^+(i)\ge |\mathcal{I}|\cdot (1-rqd/C^2)m\). 
Hence, \(\sum_{i=1}^m d_D^-(i)\ge (1-2rqd/C^2)m^2\) for sufficiently large \(C\). 
It follows by averaging that the number of indices \(i\) with \(d_D^-(i)<(1-rqd/\sqrt{C})m\) (type \ref{item:bad_cluster_T_i}) is at most \(2m/C^{3/2}\). 

We add arbitrary clusters to $\mathcal{F}_{\pi}$ to ensure $|\mathcal{F}_{\pi}|=(m-1)/C$.
By relabeling the sets $U_1',\dots,U_m'$, we may assume that $\mathcal{F}_{\pi}=\{U_i':i\in [m]\setminus [m-|\mathcal{F}_{\pi}|]\}$. 
Next, we define a bipartite graph $G_{\pi}$ between $A:=\bigcup_{U_i'\in \mathcal{F}_{\pi}} U_i'$ and $B:=\{U_2',\dots,U_m'\}\setminus \mathcal{F}_{\pi}=\{U_2',\dots,U_{m'+1}'\}$. 
Note that by the choice of \(m'\), we have 
\[|A| = (C-1)(m-1)/C=m',\quad\text{and}\quad |B|=m',\]
so both parts have size \(m'\). 
We put an edge between $v\in A$ and $U_i'\in B$ if $\delta_{\ell}(H[U_i'\cup \{v\}])\ge(\delta+2\gamma/3)\binom{|U_i'\cup \{v\}|-\ell}{k-\ell}$ and \(\mathcal{P}^{U_i'\cup \{v\}}\) is an \((F,2\beta/3,t,2c/3)\)-good partition of \(H[U_i'\cup \{v\}]\).
By $E_2$, we have $d_{G_{\pi}}(v)\ge (1-1/C^2)m-|\mathcal{F}_{\pi}|\ge 3m'/4$ for all $v\in A$. 
Since all sets satisfying \ref{item:bad_cluster_reachable} are in $\mathcal{F}_{\pi}$, for every $U_i'\in B$, we have $d_{G_{\pi}}(U_i')\ge m'-\exp{(-\theta^2 C/200)}n\ge 3m'/4 $. 
Thus, $\delta(G_{\pi})\ge 3m'/4$, and by Lemma \ref{lem:spread_PM}, there is a $(C_{\ref{lem:spread_PM}}/m')$-edge-spread distribution on perfect matchings of $G_{\pi}$.

Now we define the random sets $U_i$ as follows. First, sample $\pi$ from the uniform distribution on permutations of $[n]$ conditioning on $E_1\cap \cdots\cap E_7$. Then, sample $M_{\pi}$ from the $(C_{\ref{lem:spread_PM}}/m')$-edge-spread distribution on perfect matchings of $G_{\pi}$. For each $U_i'\in B$, let $u_i\in A$ be the vertex that it is assigned to in $M_{\pi}$, and define $U_i:=U_i'\cup \{u_i\}$. Let $U_1:=U_1'$. 
This concludes the algorithm that defines $\mathcal{U}:=\{U_1,\dots, U_{m'+1}\}$ with $m$ in the lemma statement being $m'+1$. 
Finally, we verify the properties \ref{item:(R1)}--\ref{item:(R5)} required by the theorem. 
\begin{enumerate}[label=(\arabic*)]
\item The required condition on \(|U_1|\) holds by construction. Moreover, for \(i\ge 2\), we have \(|U_i|=|U_i'|+1=C\).

\item For \(i=1\), this desired property holds since we condition on \(E_1\), while for \(i\ge 2\), it follows from the definition of the edges in the bipartite graph \(G_{\pi}\).

\item Now we verify \ref{item:(R-robust)}. 
For every \(i\in [m'+1]\), we have \(|U_i\cap V_j|\ge 2c|U_i|/3\) for all \(j\in [d]\), and the number of bad vertices in \(U_i\) is at most \(2\alpha |U_i|\ll c|U_i|\). 
In particular, this implies that \(\mathcal{T}_i\neq \emptyset\). 
For $i=1$, $d_D^-(1)\ge (1-\exp{(-\varepsilon^2C/50)})(m-1)$ by Lemma~\ref{lem:U_1} and \(E_7\); for $i\in[2,m'+1]$, we have $d_D^-(i)\ge (1-rqd/\sqrt{C})m$ since these $U_i'$ are good sets and so satisfy the negation of \ref{item:bad_cluster_T_i}.
Let $D'$ be the subgraph of $D$ induced on the first $m'+1$ vertices, and then we have $d_{D'}^-(i)\ge m'-rqdm/\sqrt{C}\ge (m'+1)/2$. Note that for every $i\in [m'+1]$, $d_{D'}^+(i)\ge m'-rqdm/C^2\ge (m'+1)/2$. 
Thus, $\delta^0(D')\ge \card{D'}/2$. By a theorem of Ghouila-Houri~\cite{ghouila1960condition}, $D'$ contains a directed Hamilton cycle, which we may assume, without loss of generality, have edges \((i,i-1\)) for all \(2\le i\le m'+1\). 
Then by the definition of $D$ and the assumption, we have for every $2\le i\le m'+1$ and every vertex $v\in T_{i}$, $\card{F^{\mu}_{\mathcal{P}}(v, U_{i-1})}\ge \card{F^{\mu}_{\mathcal{P}}(v, U_{i-1}')}\ge 2\varepsilon{\card{U_{i-1}'}}^{r-1}/3\ge \varepsilon{\card{U_{i-1}}}^{r-1}/2.$
            
\item For $i\in [m]$, \ref{item:(R4)} holds by the choice of $\mathcal{F}_{\pi}$ that contain all sets satisfying \ref{item:bad_cluster_index_verctor} and \(C\) is sufficiently large.

\item Let $y_1,\dots,y_s\in V(H)$ and $f:[s]\rightarrow[m]$. For $i\in [s]$, let $D_i$ be the event that $y_i\in U_{f(i)}$, let $D_i^1$ be the event that $y_i\in U_{f(i)}'$, and let $D_i^2$ be the event that $y_i=u_{f_i}$. Note that 
\begin{equation}\label{union bound}
    \bigcap_{i\in[s]}D_i=\bigcap_{i\in[s]}\left(D_i^1\cup D_i^2\right)=\bigcup_{S\subseteq [s]}\left(\bigcap_{i\in S}D_i^1\cap\bigcap_{i\in [s]\setminus S}D_i^2 \right).
\end{equation}
 Next, we estimate the probability that the event $\bigcap_{i\in S}D_i^1\cap\bigcap_{i\in [s]\setminus S}D_i^2$ holds. Note that $D_i^1$ is completely determined by the permutation $\pi$, and $D_i^2$ holds only if $y_i$ is matched to $U_i'$ by $M_{\pi}$. Then, for every $S\subseteq [n]$, we have
\[
\mathbb{P}\Big[\bigcap_{i\in S}D_i^1\Big]\le \Big(\frac{C^{2|S|}(n-|S|)!}{n!}\Big)/\mathbb{P}[E_1\cap\cdots \cap E_7]\le\frac{100}{93}\Big(\frac{eC^2}{n}\Big)^{|S|},
\]
            and for every $\pi'\in E_1\cap\cdots \cap E_7$, 
\[
\mathbb{P}\Big[\bigcap_{i\in [s]\setminus S}D_i^2|\pi=\pi'\Big]\le\Big(\frac{C_{\ref{lem:spread_PM}}}{m'}\Big)^{s-|S|}.
\]
Therefore, for every $S\subseteq [n]$, as $C'$ is sufficiently large, we obtain that
\[
\mathbb{P}\Big[\bigcap_{i\in S}D_i^1\cap\bigcap_{i\in [s]\setminus S}D_i^2\Big]=\mathbb{P}\Big[\bigcap_{i\in S}D_i^1\Big]\mathbb{P}\Big[\bigcap_{i\in [s]\setminus S}D_i^2\Big|\bigcap_{i\in S}D_i^1\Big]\le \frac{100}{93}\Big(\frac{eC^2}{n}\Big)^{|S|}\Big(\frac{C_{\ref{lem:spread_PM}}}{m'}\Big)^{s-|S|}\le \Big({\frac{C'}{2n}}\Big)^s.
\]
Finally, the result follows by (\ref{union bound}) and the union bound over the $2^s$ choices of $S\subseteq [s]$.
        
\end{enumerate}

The proof is complete.
\end{proof}

\section{Proof of Theorem~\ref{thm:main_technical_theorem}}\label{sec:spreadness_frm_vertex_spreadness}

In this section, we prove Theorem~\ref{thm:main_technical_theorem}.
We first introduce several additional definitions.

Recall that for a partition \(\mathcal P=\{V_1,\dots,V_d\}\) and a lattice \(L\subseteq L_{\max}^{d}\), the associated coset group is \(Q(\mathcal P,L)=L_{\max}^{d}/L.\) 
For any \(\vec{i}\in L_{\max}^{d}\), the \emph{residue} of \(\vec{i}\) in \(Q(\mathcal P,L)\) is \(R_Q(\vec{i}):=\vec{i}+L.\)
For any subset \(A\subseteq V\) with \(|A|\in r\mathbb N\), the \emph{residue} of \(A\) in \(Q\) is \(R_Q(A):=R_Q(\vec{i}_{\mathcal P}(A)).\)

Let \(q \in \mathbb{N}_0\). A (possibly empty) \(F\)-packing \(M\) in \(H\) of size at most \(q\) is a \(q\)-\emph{solution} for \((\mathcal{P}, L)\) (in \(H\)) if \(\vec{i}_{\mathcal{P}}(V(H) \setminus V(M)) \in L\); we say that \((\mathcal{P}, L)\) is \emph{\(q\)-soluble} if it has a \(q\)-solution.

The following theorem from~\cite{han2020complexity} provides a necessary and sufficient condition for \(H\) to contain a perfect \(F\)-packing, in terms of the vertex partition and robust $F$-lattice.

\begin{theorem}\cite[Theorem 3.1]{han2020complexity}\label{thm:han_structural_theorem}
         Let $k,\ell\in \mathbb{N}$ where $\ell\le k-1$ and let $F$ be an $r$-vertex $k$-graph. 
         Define $D,q,t,n_0\in\mathbb{N}$ and $\beta,\mu,\gamma,c>0$ where
         \[
         1/n_0\ll\beta,\mu\ll\gamma,c,1/r,1/D,1/q,1/t.
         \]
         Let $H$ be a $k$-graph on $n\ge n_0$ vertices where $r$ divides $n$. 
         Suppose that 
         \begin{enumerate}[label=(\roman*)]
             \item $\delta_{\ell}(H)\ge (\delta(F,\ell,D)+ \gamma)\binom{n-\ell}{k-\ell}$;\label{item:\romannumeral1}
             \item $\mathcal{P}=\{V_1, \dots, V_d\}$ is an $(F,\beta,t,c)$-good partition of $V(H)$;\label{item:\romannumeral2}
             \item $|Q(\mathcal{P},L_{\mathcal{P},F}^{\mu}(H))|\le q$.\label{item:\romannumeral3}
         \end{enumerate}
         Then $H$ contains an $F$-factor if and only if $(\mathcal{P},L^{\mu}_{\mathcal{P},F}(H))$ is $q$-soluble.
\end{theorem}

We are now ready to prove Theorem~\ref{thm:main_technical_theorem}, 
using Theorem~\ref{thm:han_structural_theorem} as a key tool.

\begin{proof}[Proof of Theorem~\ref{thm:main_technical_theorem}]
     
Suppose we have the constants satisfying the following hierarchy
\[1/n_0\ll \eta\ll 1/C'' \ll1/C'\ll 1/C\ll \beta,\mu\ll \alpha,\gamma,\varepsilon,c,1/k,1/r,1/d,1/q,1/t\le 1,\]
together with \(\alpha\ll c.\)
Let $n\ge n_0$ be a multiple of $r$, and let the $k$-graph $G$ be an $F$-packing on $n$ vertices.
Let $H$ be an $n$-vertex $k$-graph satisfying the assumptions of the theorem. 
The main task of this proof is to define $\psi:G\hookrightarrow H$, a random embedding of an $F$-factor in $H$, and we achieve it in a few steps.

\medskip
\textbf{Step 1: Sample random clusters.} Applying Lemma \ref{lem:random_cluster_degree}, we obtain a random partition $\mathcal{U}=\{U_1,\dots,U_m\}$ of $V(H)$ that satisfies properties \ref{item:(R1)}--\ref{item:(R5)}.
         
\medskip
\textbf{Step 2: Adjust the random partition.} 

Our goal is to apply Theorem \ref{thm:han_structural_theorem} to each random cluster and obtain an $F$-factor in each of them.
By \ref{item:(R-degree)}, it suffices to establish that each cluster is $q$-soluble.
For this we adjust the partition by moving at most $rq$ vertices from the vertex set $T_{i+1}$ to $U_i$ so that the solubility condition is satisfied.
Below we explain how this is achieved.
By \ref{item:(R4)}, let us first consider the lattice $L_{\mathcal{P},F}^{\mu}(H)$.

We do this in an inductive manner.
Let $i\in[m-1]$ and suppose we have corrected the first $i-1$ clusters, that is, the clusters $\{U_j':j<i\}$ satisfy $\vec{i}_{\mathcal{P}}(U_j')\in L_{\mathcal{P},F}^{\mu}(H)$ and $U_j'\subseteq U_j\cup T_{j+1}$.
Let $J_i$ denote the set of vertices removed from $U_i$ in the previous step (they were added to $U_{i-1}'$). 
Now we consider the residues of the remaining set of vertices in the $i$-th cluster with respect to $\mathcal{P}$. 
We write \(Q:=Q(\mathcal{P},L_{\mathcal{P},F}^{\mu}(H))\) for brevity. 
Note that $ \vec{i}_{\mathcal{P}}(V(H))\in L_{\mathcal{P},F}^{\mu}(H) $ and by the induction hypothesis, we have \(\vec{i}_{\mathcal{P}}(V(H)\setminus(\bigcup_{j<i}U_j'))\in L_{\mathcal{P},F}^{\mu}(H)\), namely, \(R_Q(V(H)\setminus(\bigcup_{j<i}U_j'))=\vec{0}+L_{\mathcal{P},F}^{\mu}(H)\). 
Thus, $ R_Q(U_{i}\setminus J_i)+\sum_{j> i}R_Q(U_{j})=\vec{0}+L_{\mathcal{P},F}^{\mu}(H) $. Suppose $ R_Q(U_{i}\setminus J_i)=\vec{v_0}+L_{\mathcal{P},F}^{\mu}(H) $ for some $\vec{v_0}\in L_{\max}^d$ and we get $ \sum_{j> i}R_Q(U_{j})=-\vec{v_0}+L_{\mathcal{P},F}^{\mu}(H) $.

\begin{claim}
    There exists a $q'r$-set $ J_{i+1}\subseteq  T_{i+1}$ for some $ q'\le q-1$ such that $R_Q(J_{i+1})=-\vec{v_0}+L_{\mathcal{P},F}^{\mu}(H)$.
\end{claim}
    
\begin{proof}
Since $\sum_{j> i}R_Q(U_{j})=-\vec{v_0}+L_{\mathcal{P},F}^{\mu}(H)$, by decomposing the index vectors of the clusters into sums of index vectors of $r$-sets, we may fix a family of $r$-sets $ S_1,\dots,S_{h}\subseteq\bigcup_{j> i} U_{j}$ for some $h\ge q$ such that \(\sum_{j\in[h]}R_Q(S_j)=-\vec{v_0}+L_{\mathcal{P},F}^{\mu}(H).\) 
Consider the $h+1$ partial sums $\sum_{j\in[h']}R_Q(S_j)$ for $0\le h'\le h$. 
Since $|Q(\mathcal{P},L_{\mathcal{P},F}^{\mu}(H))|\le q$, by the pigeonhole principle, there exist $j_1<j_2$ such that \(\sum_{j_1\le j \le j_2}R_Q(S_j)=\vec{0}+L_{\mathcal{P},F}^{\mu}(H).\) 
Removing such subsequences iteratively, we obtain a subfamily $ S_1,\dots,S_{q'}$ with $q'\le q-1$ satisfying \(\sum_{j\in[q']}R_Q(S_j)=-\vec{v_0}+L_{\mathcal{P},F}^{\mu}(H).\) 
Finally, since each part of $\mathcal{P}$ intersects $T_{i+1}$ in at least $rq$ vertices, we can 
choose $J_{i+1}\subseteq T_{i+1}$ of size $q'r$ such that \(R_Q(J_{i+1})=\sum_{j\in[q']}R_Q(S_j)=-\vec{v_0}+L_{\mathcal{P},F}^{\mu}(H)\).
\end{proof}

Let $U_{i}':= J_{i+1} \cup \left(U_{i}\setminus J_i\right)$, and thus $\vec{i}_{\mathcal{P}}(U_{i}')=\vec{i}_{\mathcal{P}}(J_{i+1})+\vec{i}_{\mathcal{P}}(U_{i}\setminus J_i)\in L_{\mathcal{P},F}^{\mu}(H)$. 
Finally let \(U_m':=U_{m}\setminus J_m\), and then we get \(\vec{i}_{\mathcal{P}}(U_m')\in L_{\mathcal{P}, F}^{\mu}(H)\) as $\vec{i}_{\mathcal{P}}(V(H))\in L_{\mathcal{P},F}^{\mu}(H)$. 

         \medskip
         \textbf{Step 3: Find an $F$-factor in each new cluster.} 
        Let \(J_1:=\emptyset\) and \(J_{m+1}:=\emptyset\). Then for every \(i\in [m]\), we have $U_i'=(U_i\setminus J_i)\cup J_{i+1}$, where $J_i\subseteq T_i$, $J_{i+1} \subseteq T_{i+1}$, and $|J_i|,|J_{i+1}|\leq rq$. 
        Note that by \ref{item:(R-robust)} and \(J_{m+1}=\emptyset\), for \(i\in [m]\) and every vertex $v \in J_{i+1}$, we have $\left| F^{\mu}_{\mathcal{P}}(v, U_{i}) \right| \geq \varepsilon{|U_{i}|}^{r-1}/2$.
        Thus, we find an $F$-packing $M_1^i$ in $U_i'$: for each vertex $v\in J_{i+1}$, 
        pick a copy $F_v$ of $F$ that contains $v$ and satisfies $\vec{i}_{\mathcal{P}}(V(F_v)) \in I_{\mathcal{P}, F}^{\mu}(H)$. 
        The copies of $F$ can be chosen greedily, as $H[U_i']$ contains at least $\varepsilon{|U_i'|}^{r-1}/3$ such copies of $F$ for every vertex $v \in J_{i+1}$, and we have $|V(M_1^i)| \leq r^2q$.

        Next, consider the induced subgraph $H_i:=H[U_{i}'\setminus V(M_1^i)]$ and its partition $ \mathcal{P}_i:=\mathcal{P}^{V(H_i)}$. 
        Note that $|V(H_i)|\ge |U_i|-rq-r^2q\ge |U_i|-2r^2q$ and \(|V(H_i)|\le |U_i'|\).
        We claim that $H_i$ and $ \mathcal{P}^{H_i}$ satisfy the properties described in Theorem \ref{thm:han_structural_theorem}. By applying Theorem \ref{thm:han_structural_theorem} to $ H_i$ and $ \mathcal{P}^{H_i}$, we obtain an $F$-factor $ M_2^i $ in $H_i$. 
        Then $M_1^i\cup M_2^i$ is an $F$-factor of $H[U_i']$.

        To establish the claim, we first show that $(\mathcal{P}_i,L^{\mu/3}_{\mathcal{P}_i,F}(H_i))$ is $0$-soluble, i.e.,~$\vec{i}_{ \mathcal{P}_i}(V(H_i))\in L_{ \mathcal{P}_i,F}^{\mu/3}(H_i) $. 
Indeed, since $\vec{i}_{\mathcal{P}}(U_i')\in L_{\mathcal{P},F}^{\mu}(H)$ and each copy $F_v$ in $M_1^i$ has index vector in $I_{\mathcal{P},F}^{\mu}(H)$, we have
$ \vec{i}_{\mathcal{P}}(V(H_i))=\vec{i}_{\mathcal{P}}(U_{i}'\setminus V(M_1^i))\in L_{\mathcal{P},F}^{\mu}(H) $. 
        By \ref{item:(R4)} and the definition of index vector, it follows that $\vec{i}_{ \mathcal{P}_i}(V(H_i))=\vec{i}_{\mathcal{P}}(V(H_i))\in L_{\mathcal{P},F}^{\mu}(H)\subseteq L_{ \mathcal{P}^{U_i},F}^{\mu/2}(H[U_i])$. 
        Moreover, since $|U_i\setminus V(H_i)|\le 2r^2q$, deleting these vertices removes at most $2r^2q|U_i|^{r-1}$ copies of $F$. 
        As $C$ is sufficiently large, for every $\vec{v}\in I_{ \mathcal{P}^{U_i}}^{\mu/2}(H[U_i])$, there are still at least $\frac{\mu}{2}|U_i|^r - 2r^2q|U_i|^{r-1} \ge \frac{\mu}{3}|V(H_i)|^r$ copies of $F$ in $H_i$ with index vector $\vec{v}$, which implies $\vec{v}\in I_{ \mathcal{P}_i,F}^{\mu/3}(H_i)$. Consequently, $L_{ \mathcal{P}^{U_i},F}^{\mu/2}(H[U_i])\subseteq L_{ \mathcal{P}_i,F}^{\mu/3}(H_i)$.
  
        Furthermore, since \(L_{\mathcal{P},F}^{\mu}(H)\subseteq L_{\mathcal{P}_i,F}^{\mu/3}(H_i)\) and \(|\mathcal{P}|=|\mathcal{P}_i|=d\), we have \[|Q(\mathcal{P}_i,L_{\mathcal{P}_i,F}^{\mu/3}(H_i))|=|L_{\max}^d/L_{\mathcal{P}_i,F}^{\mu/3}(H_i)|\le |L_{\max}^d/L_{\mathcal{P},F}^{\mu}(H)|=|Q(\mathcal{P},L_{\mathcal{P},F}^{\mu}(H))|\le q.\]
        In addition, by \ref{item:(R-degree)} and that $C$ is sufficiently large, we have $\delta_{\ell}(H_i)\ge \delta_{\ell}(H[U_i])-2r^2q\binom{|U_i|-\ell}{k-\ell}\ge (\delta(F,\ell,D)+\gamma/3)\binom{\card{V(H_i)}-\ell}{k-\ell} $.
        Similarly, also by \ref{item:(R-degree)}, for every $j\in[d]$, $V_j\cap V(H_i)$ is $(F,\beta/3,t)$-closed in $H_i$ and $\card{V_j\cap V(H_i)}\ge c\card{V(H_i)}/3$. 

        Therefore, for $i\in [m]$, each $H_i$ equipped with partition $\mathcal P_i$ satisfies the assumptions of Theorem \ref{thm:han_structural_theorem}, and thus has an $F$-factor.
        The union of these $F$-factors form a (random) $F$-factor of $H$.

        \medskip
        \textbf{Defining $\psi$.} 
        Now we are ready to define the random embedding function and conclude the proof.
        Fix a labeling \(V(F)=\{v_1,\dots,v_r\}\). 
Let $G$ be an $F$-factor on \([n]\) such that \(G[\{(i-1)r+1,\dots,ir\}] \cong F\) for every \(i\in[n/r]\), where the vertex \((i-1)r+j\) corresponds to \(v_j\) for each \(j\in[r]\). 
        We order the copies of \(F\) in \(M_1^1\cup M_2^1\cup M_1^2\cup M_2^2\cup\cdots\cup M_1^m\cup M_2^m\) according to the ordering of the factors. 
        For each copy \(F'\) in this ordering, fix an isomorphism \(\phi_{F'}:F\to F'\). 
        We then define an embedding \(\psi:[n]\to V(H)\) by \(\psi((i-1)r+j):=\phi_{F_i}(v_j)\) for every \(i\in[n/r]\) and \(j\in[r]\), where \(F_i\) denotes the \(i\)-th copy of \(F\) in the above ordering. 

    To complete the proof, by Proposition \ref{prop:vertex_to_edge}, it suffices to show that the random embedding $\psi$ is $(3C'/n)$-vertex-spread. Let $s\in [n]$. For every two sequences of distinct vertices $x_1,\dots,x_s\in [n]$ and $y_1,\dots,y_s\in V(H)$, we need to show that $\mathbb{P}[\psi(x_i)=y_i\text{ for }i\in[s]]\le {(3C'/n)}^s$.

We first partition $[n]$ according to the sizes of the original clusters $U_i$. 
Let $B_0:=0$ and $B_i:=\sum_{j=1}^i |U_j|$ for $i\in[m]$, and define $W_i:=[B_{i-1}+1,B_i]$. 
We claim that if $x\in W_i$, then
\[\psi(x)\in U_{i-1}\cup U_i\cup U_{i+1},\]
where we set $U_0=U_{m+1}=\emptyset$.  
To see this, let $\widehat B_0:=0$ and
$\widehat B_i:=\sum_{j=1}^i |U_j'|$, and define \(\widehat W_i:=[\widehat B_{i-1}+1,\widehat B_i].\) 
Since $M_1^i\cup M_2^i$ is an $F$-factor of $H[U_i']$, by the definition of $\psi$ we have \(\psi(\widehat W_i)=V(M_1^i\cup M_2^i)=U_i'.\)
Recall that \(U_i'=(U_i\setminus J_i)\cup J_{i+1}\) and that \(J_1=J_{m+1}=\emptyset\). 
Hence  
\[\widehat B_i
=\sum_{j=1}^i |U_j'|
=\sum_{j=1}^i (|U_j|-|J_j|+|J_{j+1}|)
=B_i+|J_{i+1}|.\]
Therefore, for every $i\in[m]$, \(W_i=[B_{i-1}+1,B_i]\subseteq \widehat W_{i-1}\cup \widehat W_i,\) where we set $\widehat W_0=\emptyset$. 
It follows that if $x\in W_i$, then 
\[\psi(x)\in U_{i-1}'\cup U_i' \subseteq U_{i-1}\cup U_i\cup U_{i+1},\] proving the claim.

Now, for each $x\in[n]$, let $w(x)$ be the unique index such that $x\in W_{w(x)}$.
Thus the event $\psi(x_i)=y_i$ for all $i\in[s]$ implies that
\[
y_i\in U_{w(x_i)-1}\cup U_{w(x_i)}\cup U_{w(x_i)+1}
\quad\text{for every }i\in[s].
\]
Consequently, by a union bound over the at most $3^s$ possible choices of clusters and by property~\ref{item:(R5)}, we obtain
\[\begin{aligned}
\mathbb{P}\left[\forall i\in[s],~\psi(x_i)=y_i\right]
&\le \mathbb{P}\left[\forall i\in[s],~y_i\in U_{w(x_i)-1}\cup U_{w(x_i)}\cup U_{w(x_i)+1}\right]\\
&\le3^s\left(\frac{C'}{n}\right)^s=\left(\frac{3C'}{n}\right)^s .
\end{aligned}\]
Thus the random embedding $\psi$ is $(3C'/n)$-vertex-spread. 

By Proposition \ref{prop:vertex_to_edge}, we conclude that there is a $(C''/n^{1/m_1(F)})$-spread distribution on subgraphs of $H$ which are isomorphic to $G$ (note that $m_1(F)=m_1(G)$).
The \emph{moreover} part of the theorem follows from Theorem \ref{thm:FKNP}.
\end{proof}

\section{Proof of Theorem~\ref{thm:main_thm_factor}}\label{sec:robustness_F}

In this subsection, we prove Theorem~\ref{thm:main_thm_factor} by applying Theorem~\ref{thm:main_technical_theorem}. 
Let $F$ be an $r$-vertex $k$-chromatic graph. 
By the definition of ${\chi}_{cr}(F)$, we have 
\begin{equation}\label{eq:chromatic_number}
    \frac{1}{{\chi}_{cr}(F)}=\frac{r-\sigma(F)}{(k-1)r}\le \frac{r-1}{(k-1)r}.
\end{equation}
To verify the assumptions required in Theorem~\ref{thm:main_technical_theorem}, we use a variant of Lemma~4.2 from \cite{lo2015f}. 
For any vertex \(u\in V(H)\), let \(W(u)\) denote the collection of \((k-1)\)-sets \(S\subseteq N(u)\) such that $S$ spans a clique in \(H\).
For any set \(T\subseteq V(H)\), define \(N(T):=\bigcap_{v\in T} N(v)\).

    \begin{lemma}\cite[Lemma 4.2]{lo2015f}\label{lem:alpha_good}
    Let $k,r\ge 2$ be integers and let $\alpha>0$. There exists a constant $\eta=\eta(k,r,\alpha)$ such that the following holds for sufficiently large $n$. 
    Let \(F\) be a \(k\)-chromatic graph on \(r\) vertices.
    For every \(n\)-vertex graph \(G\), two vertices \(u,v\in V(G)\) are $(F,\eta,1)$-reachable to each other if the number of \((k-1)\)-sets \(S\in W(u)\cap W(v)\) with \(\card{N(S)}\ge \alpha n\) is at least \(\alpha \binom{n}{k-1}\).
    \end{lemma}

The following lemma bounds the size of the associated coset group. 

\begin{lemma}\cite[Proposition 9.3]{han2020complexity}\label{lem:coset_group_of_factor}
Let \(k,r,t,n_0\in \mathbb{N}\) where \( k \geq 2\) and let \(\beta,\mu,\gamma>0 \) so that
\[
1/n_0\ll \beta, \mu\ll \gamma\ll 1/r, 1/t.
\]
Let \(F\) be an unbalanced \(r\)-vertex \(k\)-chromatic graph. Suppose \(G\) is a graph on \( n \geq n_0 \) vertices such that \( \delta(G) \geq (1 -1/\chi_{cr}(F)+\gamma)n \) with an \( (F, \beta, t, 1/r) \)-good partition \( \mathcal{P} \) where \(|\mathcal{P}|=d\). Then \(|Q(\mathcal{P}, L_{\mathcal{P}, F}^{\mu}(G))|\le (2r-1)^d.\)
\end{lemma}

    To have a more precise threshold for strictly 1-balanced $F$, we consider the spread properties of the \(F\)-graph of \(G\). 
    For hypergraphs \(F\) and \(G\), the \(F\)-\emph{graph} of \(G\), denoted by \(G_F \), is the \(\card{V(F)}\)-uniform multi-hypergraph with vertex set \(V(G)\) in which every copy of \(F\) in \(G\) corresponds to a distinct hyperedge of \(G_F\) on the same set of vertices. 
    Note that \(G\) has an \(F\)-factor if and only if \(G_F\) contains a perfect matching. 
    Let \( \mathbb{G}_F(n,p) \) denote the binomial random multi-hypergraph on \(n\) vertices where every edge of the \(F\)-graph of the complete hypergraph is included independently with probability \(p\). 


    \begin{theorem}\cite[Theorem 18]{riordan2022random}\label{couple}
        For every \(r\in \mathbb{N}\), there exists \( a=a_{\ref{couple}}(r) \) such that the following holds.
        Let \(F\) be a fixed strictly 1-balanced graph. 
        If \( p = p(n) \leq \log^{2}(n)/n^{1/m_1(F)} \), then, for some \( \pi = \pi(n) \sim ap^{e_F }\), there exists a coupling for \(G=\mathbb{G}\left(n,p\right)\) and \(G_F=\mathbb{G}_{F}(n,\pi)\) such that, a.a.s.~for every \(F\)-edge present in \(G_F\) the corresponding copy of F is present in G.
    \end{theorem}

Now we are ready to prove Theorem~\ref{thm:main_thm_factor}.

\begin{proof}[Proof of Theorem~\ref{thm:main_thm_factor}]
Define additional constants \(\alpha,\varepsilon,\theta>0\) such that 
\[ 1/n\ll\beta, \mu\ll \alpha,\varepsilon\ll \theta,\gamma\ll1/r,1/k. \]
To prove the theorem we consider two cases. 
         
We first prove the last part of the theorem -- if \(G\) contains an \(F\)-factor then $G$ contains an \(F\)-packing $M_0$ of size at most $(2r-1)^r$ such that $ \vec{i}_{\mathcal{P}}(V(G)\setminus V(M_0))\in L_{\mathcal{P},F}^{\mu}(G)$. 
Note that \(\card{\mathcal{P}}\le r\), by Lemma~\ref{lem:coset_group_of_factor}, we have $\card{Q(\mathcal{P},L_{\mathcal{P},F}^{\mu}(G))}\le (2r-1)^r.$ 
It then follows from Theorem~\ref{thm:han_structural_theorem} that $(\mathcal{P},L^{\mu}_{\mathcal{P},F}(G))$ is $(2r-1)^r$-soluble, which implies the existence of a desired $F$-packing.

Conversely, assume that there exists an \(F\)-packing \(M_0\) of size at most \((2r-1)^r\) satisfying $ \vec{i}_{\mathcal{P}}(V(G)\setminus V(M_0))\in L_{\mathcal{P},F}^{\mu}(G)$. 
Let \(G':= G-V(M_0)\), and its partition \(\mathcal{P}' := \mathcal{P}^{G'}\). Set \(m := |V(G')|\).
We aim to apply Theorem~\ref{thm:main_technical_theorem} to \(G'\) and \(\mathcal{P}'\) with \(\ell=1\) and \(D=5r^2\).
Recall that $\delta(F,1,5r^2)=1-1/{\chi}_{cr}(F)$ by Theorem \ref{thm:deltaFD}.
Since removing \(V(M_0)\) deletes at most $r(2r-1)^r$ vertices, the minimum degree and the number of reachable sets between any two vertices are only slightly reduced. 
Hence, \(\delta(G')\ge (1-1/{\chi}_{cr}(F)+\gamma/2)m\), \(\mathcal{P}'\) is an \((F,\beta/2,2^{h-1},1/r)\)-good partition, and Lemma~\ref{lem:coset_group_of_factor} implies that $|Q(\mathcal{P}',L_{\mathcal{P}',F}^{\mu}(G'))|\le (2r-1)^r$. 

It remains to verify assumption (iii) of Theorem~\ref{thm:main_technical_theorem}. 
We show that for every vertex \(v\in V(G')\), there exist at least \(\theta m^{r-1}\) \((r-1)\)-sets \(R \subseteq V(G')\) such that \(R\cup \{v\}\) spans a copy of \(F\) in \(G'\). 
Indeed, since \(\delta(G')\ge (1-1/{\chi}_{cr}(F)+\gamma/2)m\), it follows from \eqref{eq:chromatic_number} that for every \((k-1)\)-set \(S\), \(|N(S)|\ge (1/r+(k-1)\gamma/2)m\). 
Consequently, every vertex \(v\) satisfies \(|W(v)|\ge \frac{1}{r^{k-1}}\binom{m}{k-1}\). 
Fix \(v \in V(G')\), and let \(G''\) be the graph obtained from \(G'\) by adding a duplicate vertex \(v'\) with \(N(v')=N(v)\).
By Lemma~\ref{lem:alpha_good}, the vertices \(v\) and \(v'\) are \((F,\theta,1)\)-reachable in \(G''\). 
Hence there exist at least \(\theta m^{r-1}\) \((r-1)\)-sets \(R \subseteq V(G')\) such that both \(R\cup\{v\}\) and \(R\cup\{v'\}\) span copies of \(F\) in \(G''\). 
Since \(v'\) is a duplicate of \(v\), each such set \(R\) yields a copy of \(F\) containing \(v\) in \(G'\).

Let $ B_0 $ be the set of all vertices $v$ such that \(\card{F_{\mathcal{P}'}^{\mu}(v)}\le \varepsilon m^{r-1}\).
Note that for each \(v\in B_0\), the vertex \(v\) lies in at least \((\theta-\varepsilon)m^{r-1}\ge \theta m^{r-1}/2\) copies of \(F\) whose index vector is not in \(I_{\mathcal{P}',F}^{\mu}(G')\), yielding
\[|B_0|\le\frac{r\binom{r+d-1}{r}\mu m^r}{\theta m^{r-1}/2}\le (2r)^{r}\mu m/\theta <\alpha m,\]
where the last inequality follows from $\mu\ll\theta,\alpha, 1/r$. 
Therefore, assumption (iii) of Theorem~\ref{thm:main_technical_theorem} holds.
We apply Theorem \ref{thm:main_technical_theorem} to $G'$ with $\gamma/2, \beta/2$ in place of $\gamma,\beta$, and obtain a $\left(C_1/n^{1/m_1(F)}\right)$-spread distribution on the set of $F$-factors in $G'$, which implies that a.a.s.~$G_p'$ contains an $F$-factor if \(p\ge Cn^{-1/m_1(F)}\log n\).

We now turn to the case when \(F\) is strictly 1-balanced. It suffices to prove the result for \(p=Cn^{-1/m_1(F)}\log^{1/e_F}n\).  
Denote the \(F\)-multi-hypergraph of \(G'\) by \(\mathcal{G}'\).
As there is an $\left(C_1/n^{1/m_1(F)}\right)$-spread distribution $\lambda$ on the family of $F$-factors in $G'$, it gives rise to a \(\left(C_1^{e_F}/n^{r-1}\right)\)-spread distribution $\lambda'$ on the perfect matchings of \(\mathcal{G}'\). 
Indeed, suppose $S$ is a matching of size $s$ in \(\mathcal{G}'\) (if $S$ does not form a matching then $\lambda'(S)=0$), then we have
\[
\lambda'(S) = \lambda(T)\le C_1^{e_F\cdot s} n^{-\frac{1}{m_1(F)}e_F\cdot s} = C_1^{e_F\cdot s} n^{-(r-1)s},
\]
where $T$ is the $F$-packing of size $s$ in $G'$ corresponding to $S$ and we used $|T|=e_F\cdot s$ and $m_1(F)=e_F/(r-1)$.
Therefore, Theorem \ref{thm:FKNP} implies that a.a.s.~\(\mathcal{G}'_{\pi}\) contains a perfect matching if \(\pi\ge K_{\ref{thm:FKNP}}C_1^{e_F}\log n/ n^{(r-1)}=ap^{e_F}\). Thus, by Theorem \ref{couple}, a.a.s.~\(G_p'\) contains an \(F\)-factor.
\end{proof}

\section{Separating a sub-lattice}\label{sec:lattice_preparation}

In this section we include a simple fact from the lattice theory which will be used to separate two families of robust index vectors later in the proof. 

\begin{fact}\label{fact:separating_sublattice}
Suppose $I_1, I_2$ are finite sets of vectors in $\mathbb Z^d$.
Let $L_1$, $L_2$ and $L$ be lattices generated by $I_1$, $I_2$ and $I_1\cup I_2$, respectively.
Then for every $\vec{x}\in L$, there exists $\vec{y}\in L_2$ such that $\vec{x} - \vec{y}\in L_1$.
\end{fact}

\begin{proof}
    
    For $\vec{x}\in L$, as $L$ is generated by $I_1\cup I_2$, we have the integer equation 
    \[
    \vec{x}=\sum_{\vec{v}\in I_1\cup I_2}c_{\vec{v}}\vec{v},
    \]
    where $c_{\vec{v}}\in \mathbb{Z}$, $\vec{v}\in I_1\cup I_2$ are the integer coefficients.
    Then let $\vec{y}=\sum_{\vec{v}\in I_2}c_{\vec{v}}\vec{v}\in L_2$ and note that $\vec{x}- \vec{y} = \sum_{\vec{v}\in I_1}c_{\vec{v}}\vec{v} \in L_1$, as desired.
    
\end{proof}

In the applications below, \(\vec{x}\) will be the index vector of the current remaining vertex set, while \(\vec{y}\) will be the index vector of a first-stage matching whose edges have index vectors from the second family. 
Thus the task is not merely to choose edges of the second type, but to choose their total index vector so that it removes the second-lattice contribution of \(\vec{x}\). 
This is achievable because we know the (finite) set of generators for $L_2$.
Once this is achieved, the remaining vertex set has index vector in \(L_1\). 
The main difficulty is to realize this prescribed lattice correction by an actual matching while preserving the spread estimates needed later.
This is the key new idea in the proof of Theorem \ref{thm:auxiliary_thm_perfect_matching}.

\section{Proof of Theorem~\ref{thm:algorithm_PM}}\label{sec:proof_main_result}

In this section, we prove Theorem~\ref{thm:algorithm_PM}.
We first recall the structural framework developed by Gan and Han~\cite{gan2025keevash}, which reduces the existence of perfect matchings to a lattice-based solubility condition.
After introducing the necessary structural tools, we present our algorithm and prove its correctness.

\subsection{Structural tools}\label{sec:outline}

Under minimum \(\ell\)-degree conditions, the resulting partition may contain certain exceptional parts that do not satisfy strong closure properties. 
To handle this issue, Gan and Han~\cite{gan2025keevash} introduced a refined notion of robust edge-vectors, which distinguishes the roles of different parts of the partition.
We say that \(\vec{v}\in\mathbb Z^d\) is a \(k\)-vector if all coordinates of \(\vec{v}\) are non-negative and \(|\vec{v}|=k\). 

\begin{definition}[\(\mu\)-robust vectors,~\cite{gan2025keevash}]\label{def:gan_han_def}
    For \(d>s\ge 0\), let \(\mathcal{P} = \{V_0, V_1, \ldots, V_s, V_{s+1}, \ldots, V_d\}\) be a partition of \(V(H)\). 
    Given \(\mu > 0\), define \(I_{\mathcal{P}}^{s,\mu}(H) := I_{\mathcal{P},1}^{s,\mu}(H) \cup I_{\mathcal{P},2}^{s,\mu}(H)\) as the union of the following two sets: 
    \begin{enumerate}
    \item[(1)] the set \(I_{\mathcal{P},1}^{s,\mu}(H)\) consists of all \(k\)-vectors \(\vec{v} \in \mathbb{Z}^d\) such that \(\vec{v}|_i = 0\) for \(i \in [s]\), \(\sum_{j=s+1}^d \vec{v}|_j = k\) and \(H\) contains at least \(\mu n^k\) edges \(e\) with \(\vec{i}_{\mathcal{P}}(e) = \vec{v}\);
    \item[(2)] the set \(I_{\mathcal{P},2}^{s,\mu}(H)\) consists of all \(k\)-vectors \(\vec{v} \in \mathbb{Z}^d\) such that \(\vec{v}|_i = 1\) for exactly one \(i \in [s]\), \(\sum_{j=s+1}^d \vec{v}|_j = k - 1\) and every vertex \(v \in V_i\) is in at least \(\mu n^{k-1}\) edges \(e\) with \(\vec{i}_\mathcal{P}(e) = \vec{v}\).
\end{enumerate}
\end{definition}

Compared with Definition~\ref{def_lattice}, the new feature is the family
\(I^{s,\mu}_{\mathcal P,2}(H)\), which accounts for edges using one vertex from one of the parts \(V_1,\dots,V_s\). 
This refinement is crucial under general minimum \(\ell\)-degree conditions. 

Gan and Han~\cite{gan2025keevash} established the following partition lemma, which produces a refined partition tailored for minimum $\ell$-degree settings.

\begin{lemma}\cite[Lemma 2.2]{gan2025keevash}\label{lem:partition_l_degree}
Let \(k,\ell,h,m\in\mathbb N\) with \(1\le \ell\le k-1\),
and let \(\delta>0\). Suppose $1/n_0 \ll \mu \ll \beta \ll \delta' \ll \delta, 1/k, 1/m$. 
Given an $n$-vertex $k$-graph $H$ with $n \ge n_0$ and $\delta_\ell(H) \ge \delta \binom{n-\ell}{k-\ell}$, there is a partition $\mathcal{P}$ of $V(H)$ as
\[
\mathcal{P} = \{V_0, V_1, \dots, V_s, V_{s+1}, \dots, V_d\}
\]
such that with $c := \lfloor 1/\delta \rfloor$:
\begin{enumerate}
    \item $s \le 2^{\binom{c+k-2}{k-1}}$ and $d - s \le c$;
    \item $|V_0| \le k^{2^{\binom{c+k-2}{k-1}}} \left( kh + \binom{c+k-2}{k-1}m\right)$ and $\bigl|\bigcup_{0 \le i \le s} V_i\bigr| \le c\delta'n$;
    \item for $1 \le i \le s$, $|V_i| \ge (k-1)|V_0| + kh+  \binom{c+k-2}{k-1}m$;
    \item for $1 \le i \le s$, there exists $\vec{v} \in I_{\mathcal{P},2}^{s,\mu}(H)$ such that $\vec{v}|_{i} = 1$;
    \item for $s+1 \le i \le d$, $|V_i| \ge \delta'n/2$ and $V_i$ is $(\beta, 2^c)$-closed in $H\bigl[\bigcup_{s+1 \le i \le d
    } V_i\bigr]$.
\end{enumerate}
In particular, such a partition $\mathcal{P}$ of $H$ can be found in time $O(n^{2^ck+1})$.
\end{lemma} 


The conclusions of Lemma~\ref{lem:partition_l_degree} play different roles in our argument.
Item~(1) gives bounds on the numbers of exceptional and non-exceptional parts in the partition. 
Item~(2) shows that the exceptional set \(V_0\) has constant size, which allows the algorithm to search for a matching covering \(V_0\) by brute force. 
Item~(3) ensures that the sets \(V_1,\dots,V_s\) are sufficiently large, providing enough room for handling vectors in \(I_{\mathcal P,2}^{s,\mu}(H)\) in later arguments.
Item~(4) guarantees that vertices in each \(V_i\), \(i\in [s]\), can be matched in many different ways.
Finally, item~(5) provides the closure properties needed for the absorption argument. 
Moreover, the proof of Lemma~\ref{lem:partition_l_degree} in~\cite{gan2025keevash} yields a useful converse to item~(4), which we record for later use: 

\begin{enumerate}[label=$(\star)$]
    \item For every \(1 \le i \le s\) and every \(k\)-vector \(\vec{v}\notin I_{\mathcal P,2}^{s,\mu}(H)\) satisfying \(\vec{v}|_i = 1\) for exactly one \(i \in [s]\) and \(\sum_{j=s+1}^d \vec{v}|_j = k - 1\), every vertex \(v\in V_i\) is contained in fewer than \(\mu n^{k-1}\) edges \(e\) with \(\vec{i}_{\mathcal P}(e)=\vec{v}\).\label{item:extra_propert}
\end{enumerate}
Although property~\ref{item:extra_propert} is not stated explicitly
in~\cite{gan2025keevash}, it is an immediate consequence of the way the sets \(V_1,\dots,V_s\) are constructed in the proof of Lemma~\ref{lem:partition_l_degree}.

Building on the partition lemma, Gan and Han~\cite{gan2025keevash} established the following structural theorem. 
We define \(L^{s,\mu}_{\mathcal{P}}(H)\) to be the lattice in \(\mathbb{Z}^d\) generated by \(I_{\mathcal{P}}^{s,\mu}(H)\). 
We use the analogous notion of \(q\)-solution for matchings. 
For later use, we strengthen it as follows.
Given a set \( U \subseteq V(H) \), we define that \( (\mathcal{P}, L) \) is \emph{\((U, q)\)-soluble} if there is a matching \( M \) in \( H \) such that \( M \) covers \( U \) and \( M \) is a \( (|U| + q) \)-solution.

Note that both the theorem below and Theorem \ref{thm:auxiliary_thm_perfect_matching} output a parameter $m$, which we use as an input to Lemma \ref{lem:partition_l_degree}.

\begin{theorem}\cite[Theorem 2.3]{gan2025keevash}\label{thm:gan_han_structural}
Let $k, \ell, q \in \mathbb{N}$ where $\ell \le k-1$ and let $\gamma > 0$ be given. 
There exist $n_0, m := m_{\ref{thm:gan_han_structural}}(k, q) \in \mathbb{N}$ and $\delta', \beta, \mu > 0$ such that
\[
1/n_0 \ll \beta, \mu \ll \delta' \ll \gamma, c_{k,\ell}^*, 1/k,1/q.
\]
Let $H$ be an $n$-vertex $k$-graph with \(\delta_\ell(H) \ge (c_{k,\ell}^* + \gamma) \binom{n-\ell}{k-\ell},\)
where $n \ge n_0$ and $k$ divides $n$. 
Suppose $\mathcal{P}$ is a partition of $V(H)$ satisfying Lemma~\ref{lem:partition_l_degree} (1)-(5) with $\delta = c_{k,\ell}^*$, $m$ and \(h=q\). 
Moreover, suppose $|Q(\mathcal{P}, L_{\mathcal{P}}^{s,\mu}(H))| \le q$. 
Then $H$ contains a perfect matching if and only if $(\mathcal{P}, L_{\mathcal{P}}^{s,\mu}(H))$ is $(V_0, q)$-soluble. 
Furthermore, there is an algorithm with running time $O(n^{2^{k+2}}k^2)$ that finds a perfect matching in $H$ if $(\mathcal{P}, L_{\mathcal{P}}^{s,\mu}(H))$ is $(V_0, q)$-soluble.
\end{theorem}

To apply Theorem~\ref{thm:gan_han_structural}, we also need a bound on the size of the relevant coset group. 
The following lemma, adapted from Proposition~4.1 of~\cite{gan2025keevash}, provides such a bound.
Although the original statement assumes all conditions of Lemma~\ref{lem:partition_l_degree}, its proof only uses (2) and (5), and moreover gives the additional bound $|L_{\max}^{d-s}/L_{\mathcal{P},1}^{s,\mu}(H)|\le (2k+1)^{d-s}$, which will be used later in the proof of the auxiliary theorem.

\begin{lemma}\cite[Proposition~4.1]{gan2025keevash}\label{lem:coset_group_of_l_degree}
    Suppose \(1/n_0\ll\mu \ll\delta'\ll\gamma\ll 1/k.\) Let \(H\) be an \(n\)-vertex \(k\)-graph with \(\delta_{\ell}(H)\ge (c^*_{k,\ell}+\gamma)\binom{n-\ell}{k-\ell}\), where \(n\ge n_0\) and \(k\mid n\). Suppose that \(\mathcal{P}=\{V_0,V_1,\dots,V_s,V_{s+1},\dots,V_d\}\) is a partition of \(V(H)\) satisfying Lemma~\ref{lem:partition_l_degree} (2) and (5) with \(\delta=c^*_{k,\ell}\). Then  $|Q(\mathcal{P},L_{\mathcal{P}}^{s,\mu}(H))|\le (2k+1)^{d-s}$ and $|L_{\max}^{d-s}/L_{\mathcal{P},1}^{s,\mu}(H)|\le (2k+1)^{d-s}$.
\end{lemma}

The main new ingredient of the present paper is the following auxiliary theorem. It shows that once a suitable solution is found in the random subhypergraph, a.a.s. the remaining part of the hypergraph still contains a perfect matching. 
This theorem is the key tool that allows us to transfer the deterministic lattice-based framework of Gan and Han to the random sparsification setting. 

\begin{theorem}\label{thm:auxiliary_thm_perfect_matching}
Let $k, \ell\in \mathbb{N}$ where $\ell \le k-1$ and let $\gamma > 0$ be given. 
There exist $C,n_0, m := m_{\ref{thm:auxiliary_thm_perfect_matching}}(k) \in \mathbb{N}$ and $\eta, \delta', \beta, \mu > 0$ such that
\[
1/n_0 \ll \eta\ll 1/C\ll\beta, \mu \ll \delta' \ll \gamma, c_{k,\ell}^*,1/k.
\]
Let $H$ be an $n$-vertex $k$-graph with $n \ge n_0$ and $k$ divides $n$. Suppose that 
\begin{enumerate}[label=$(S\arabic*)$]
    \item For every \(v\in V(H)\), all but at most \(\eta\binom{n-1}{\ell-1}\) sets \(S\in\binom{V(H)\setminus\{v\}}{\ell-1}\) satisfy \(\deg_H(S\cup\{v\})\ge(c_{k,\ell}^*+\gamma)\binom{n-\ell}{k-\ell}\);\label{item:(C*_degree)}
    \item $\mathcal{P}$ is a partition of $V(H)$ satisfying Lemma~\ref{lem:partition_l_degree} (1)-(5) with $\delta = c_{k,\ell}^*$, \(V_0=\emptyset\), and \(h=0\);\label{item:(C*_partition)}
    \item \(H\) and \(\mathcal P\) satisfy~\ref{item:extra_propert} with upper bound \(4\mu n^{k-1}\) in place of \(\mu n^{k-1}\). \label{item:star}
\end{enumerate}
If $\vec{i}_{ \mathcal{P}}(V(H))\in L_{ \mathcal{P}}^{s,\mu}(H)$,  then there exists an $O(1/n^{k-1})$-spread distribution on the set of perfect matchings in $H$. 
In particular, if \(p\ge C\log n/n^{k-1}\), then a.a.s.\ \(H_p\) contains a perfect matching.
\end{theorem}

\subsection{Proof of Theorem \ref{thm:algorithm_PM}}

In this subsection, we prove Theorem~\ref{thm:algorithm_PM}. 
First we present our algorithm for Theorem \ref{thm:algorithm_PM} (see Procedure \ref{alg:perfect_matching}). 
For convenience, write \(I':= I_{\max}^d \setminus I_{\mathcal{P}}^{s,\mu}(H).\)

\begin{procedure}[ht]
\caption{PerfectMatching()}
\label{alg:perfect_matching}
\SetKwInOut{Input}{Input}
\SetKwInOut{Output}{Output}

\Input{An $n$-vertex $k$-graph $H$ with $\delta_{\ell}(H)\ge (c^*_{k,\ell}+\gamma)\binom{n-\ell}{k-\ell}$ and $k\mid n$, and random subgraph $H_p$.}
\Output{a.a.s.~$H_p$ has a perfect matching or none at all.}

Choose constants $1/n_0 \ll \eta \ll 1/C \ll \beta, \mu \ll \delta' \ll \gamma , c_{k,\ell}^*, 1/k$, and set $q:=(2k+1)^k$\;

\If{$n<n_0$}{
    Check whether $H_p$ has a perfect matching by brute force, return result and halt\;
}

Apply Lemma \ref{lem:partition_l_degree} with $\delta = c_{k,\ell}^*$ and \(h=q\) to find a partition $\mathcal{P}$ of $V(H)$ and $L^{s,\mu}_{\mathcal{P}}(H)$\;

\For{every $\vec{v}\in I'$}{
    compute $V_{\vec{v}}^+$ and expose edges to obtain $E(\vec{v})$\;
}

\For{every minimal matching $M_0$ containing $V_0$ and $M_0 \subseteq H_p$}{
    \For{
    every choice of integers $(c_{\vec{v}})_{\vec{v} \in I'}$ with $c_{\vec{v}} \ge 0$ such that $\sum_{\vec{v}} c_{\vec{v}} \le q$ and $\vec{i}_{\mathcal{P}}(V(H)\setminus V(M_0)) - \sum_{\vec{v} \in I'} c_{\vec{v}} \vec{v} \in L_{\mathcal{P}}^{s,\mu}(H)$}{
        
    \For{every choice of integers $c_{\vec{v}}^* \le c_{\vec{v}} $ and edge sets $M_{\vec{v}}  \subseteq E(\vec{v})$ with $|M_{\vec{v}}|=c_{\vec{v}}^*$}{Let $M_1 := \bigcup_{\vec{v}} M_{\vec{v}} $\;
    
    \If{$M_0 \cup M_1$ is a matching}{
        \If{there exist disjoint sets $S_{\vec{v}}  \subseteq V_{\vec{v}}^+\setminus V(M_0 \cup M_1)$ with $|S_{\vec{v}}| = c_{\vec{v}} - c_{\vec{v}}^*$ for all $\vec{v} \in I'$}{
            \Return{``a.a.s.~$H_p$ contains a perfect matching'' and halt}\;
}
}
}
}
}

\Return{``$H_p$ contains no perfect matching'' and halt}\;
\end{procedure} 


\begin{proof}[Proof of Theorem~\ref{thm:algorithm_PM}]
We show correctness of Procedure \ref{alg:perfect_matching} and estimate its running time.  

Note that $c_{k,\ell}^* \ge c_{k,k-1}^* = 1/k$, so that \(c = \left\lfloor 1/c_{k,\ell}^* \right\rfloor \le k.\)
Set \(h=q := (2k+1)^k\) and $m := \max\{m_{\ref{thm:gan_han_structural}}, m_{\ref{thm:auxiliary_thm_perfect_matching}}\}$. 
We apply Lemma~\ref{lem:partition_l_degree} with this choice of $h$ and $m$, which finds a partition $\mathcal{P}$ satisfying properties (1)--(5) in time $O(n^{2^k+1})$. 
By Lemma~\ref{lem:coset_group_of_l_degree}, \(|Q(\mathcal{P}, L_{\mathcal{P}}^{s,\mu}(H))| \le (2k+1)^{d-s} \le q.\)


By Theorem \ref{thm:gan_han_structural}, if $H_p$, as a subgraph of $H$, does not contain a matching $M$ covering $V_0$ which is a $(|V_0|+q)$-solution, then $H_p$ contains no perfect matching.  
The algorithm therefore first searches for such a solution. 
If no solution exists, the algorithm outputs correctly that no perfect matching exists. Otherwise, assume a $(|V_0|+q)$-solution is found. We show that a.a.s.\ it can be extended to a perfect matching using Theorem~\ref{thm:auxiliary_thm_perfect_matching}.

\medskip
\textbf{Finding a $(|V_0|+q)$-solution.}
We first perform a preprocessing step (line 5).  
For every $k$-vector $\vec{v} \notin I_{\mathcal{P}}^{s,\mu}(H)$, we examine the edges of $H$ and define
\[
V_{\vec{v}}^+ := \Bigl\{ v \in V(H) : 
\bigl|\{ e \in E(H) : v \in e,\ \vec{i}_{\mathcal P}(e)=\vec{v} \}\bigr|
\ge \eta n^{k-1} \Bigr\}.
\]
This can be done by screening the edges of $H$ in $O(n^k)$ time.
For every vertex $v \in V(H)\setminus V_{\vec{v}}^+$, we expose all edges in $H_p$ containing $v$ with index vector $\vec{v}$, except those intersecting \(V_{\vec{v}}^+\).
Let $E(\vec{v})$ denote the collection of all such exposed edges that survive in $H_p$.  
By construction, for each \(\vec v\in I'\), we expose edges of type \(\vec v\) only at vertices outside \(V^+_{\vec v}\). Hence each vertex outside \(V_0\) is incident to at most \(\binom{d+k-1}{k}\eta n^{k-1}+|V_0|n^{k-2}\) revealed edges.

We now search for a $(|V_0|+q)$-solution in $H_p$, and the idea is to exhaust all possibilities of (small) matchings using edges either intersecting $V_0$ or from $E(\vec{v})$. 
First, we expose all edges of \(H_p\) that are incident to \(V_0\), and enumerate all minimal matchings \(M_0\) among these exposed edges that survive in \(H_p\) and cover \(V_0\). 
This can be done in time \(O(n^{(k-1)|V_0|})\), and by minimality every edge of such an \(M_0\) is incident to \(V_0\).  
Fix such a matching $M_0$. We extend it to a larger matching using edges from $\bigcup_{\vec{v} \in I'} E(\vec{v})$.
To this end, we iterate over all choices of nonnegative integers $c_{\vec{v}}$ for each $\vec{v} \in I'$,
such that $\sum_{\vec{v}} c_{\vec{v}} \le q$ and  
\[
\vec{i}_{\mathcal{P}}(V(H)\setminus V(M_0)) - \sum_{\vec{v} \in I'} c_{\vec{v}} \vec{v} \in L_{\mathcal{P}}^{s,\mu}(H).
\]
These vectors describe the types of edges that may be added to extend \(M_0\) to a valid solution. 
Fix such an assignment $(c_{\vec{v}})_{\vec{v} \in I'}$.  
For every choice of integers $c_{\vec{v}}^* \le c_{\vec{v}} $ and edge sets $M_{\vec{v}}  \subseteq E(\vec{v})$ with $|M_{\vec{v}}|=c_{\vec{v}}^*$, we define
\[
M_1 := \bigcup_{\vec{v}\in I'} M_{\vec{v}}.
\]
If $M_0 \cup M_1$ is a matching,  then \(\vec{i}_{\mathcal{P}}(V(M_1)) = \sum_{{\vec{v}}} c_{\vec{v}}^* \vec{v},\) and we attempt to extend it further (line 10).
Given such a matching $M_1$, we seek sets of vertices
\begin{equation}
\label{eq:Sv}
S_{\vec{v}} \subseteq V_{\vec{v}}^+ \setminus V(M_0 \cup M_1), \quad |S_{\vec{v}}| = c_{\vec{v}}-c_{\vec{v}}^*,
 \end{equation}
for $\vec{v}\in I'$
such that these sets are pairwise disjoint. 
If such a collection exists (line 11), then we can a.a.s. extend $M_0 \cup M_1$ to a $(|V_0|+q)$-solution using the following claim. 
In the following claim, we condition on all information exposed before line~12, as well as on the fixed choice of the sets \((S_{\vec v})_{\vec v\in I'}\) satisfying~\eqref{eq:Sv}. 
\begin{claim}\label{claim:vec_v}
With high probability, for every \(\vec{v}\in I'\) and \(v\in S_{\vec{v}}\), there exist at least
\((|V_0|+q)k\) edges of \(H_p\) containing \(v\) and having index vector
\(\vec{v}\), which are pairwise disjoint outside \(v\).
\end{claim}

\begin{proof}
Fix \(\vec{v}\in I'\) and \(v\in S_{\vec{v}}\). 
Since \(S_{\vec{v}} \subseteq V_{\vec{v}}^+\), there are at least \(\eta n^{k-1}\) \((k-1)\)-sets \(S \subseteq V(H)\setminus\{v\}\) such that \(S \cup \{v\} \in E(H)\) has index vector \(\vec{v}\). 

Let \(h:=(|V_0|+q)k\), and partition \(V(H)\setminus\{v\}\) uniformly at random into \(h\) parts \(U_1,\dots,U_h\) of nearly equal size. 
For each \(i\in[h]\), let \(X_i\) be the number of \((k-1)\)-sets \(S \subseteq U_i\) such that \(S \cup \{v\} \in E(H)\) has index vector \(\vec{v}\), and let \(Y_i\) be the number of such sets \(S\) with \(S \cup \{v\} \in E(H_p)\). 
By Lemma \ref{lem:application_McDiarmid} applied to each \(U_i\), we have 
\[
\mathbb{P}\left[X_i<\frac{\eta}{2} \card{U_i}^{k-1}\right]\le \exp{\left(-\frac{{\eta}^2n}{10h}\right)}=o(1).
\] 
Conditioning on \(X_i\ge \frac{\eta}{2} \card{U_i}^{k-1}\), we have \(\mathbb{E}[Y_i]\ge \frac{\eta }{2}\card{U_i}^{k-1}p\), so by Chernoff's bound, \(\mathbb{P}[Y_i=0]\le \exp{\left(-\frac{\eta\card{U_i}^{k-1}p}{4}\right)}=o(1).\) 
Thus \(\mathbb{P}[Y_i=0]=o(1)\), and by a union bound over all \(i\in[h]\), with probability \(1-o(1)\), for every \(i\in [h]\) there exists a set \(S \subseteq U_i\) such that \(S \cup \{v\} \in E(H_p)\). 
These yield \(h=(|V_0|+q)k\) edges in \(H_p\) containing \(v\), which are pairwise disjoint outside \(v\) by construction.

Finally, since the total number of pairs \((\vec{v},v)\) is at most \(\sum_{\vec{v}} |S_{\vec{v}}|\le q\), a union bound implies that the statement holds simultaneously for all \(\vec{v}\in I'\) and \(v\in S_{\vec{v}}\) with high probability.
\end{proof}

If for a given choice of coefficients $(c_{\vec{v}})$ and all choices of
$c_{\vec{v}}^* \le c_{\vec{v}}$ together with edge sets
$M_{\vec{v}} \subseteq E(\vec{v})$,
there do not exist pairwise disjoint vertex sets
$(S_{\vec{v}})_{\vec{v} \in I'}$
with $|S_{\vec{v}}| = c_{\vec{v}} - c_{\vec{v}}^*$,
then this choice of coefficients does not yield a valid extension,
and we proceed to the next one. 
If all choices of $M_0$ and all choices of coefficients $(c_{\vec{v}})$ are exhausted and no valid extension is found, then $H_p$ contains no $(|V_0|+q)$-solution. 
Therefore, $H_p$ contains no perfect matching by Theorem~\ref{thm:gan_han_structural}. 
In this case Procedure \ref{alg:perfect_matching} halts at line 14 with correct output.

Now suppose the algorithm finds a family of sets $(S_{\vec{v}})_{\vec{v} \in I'}$ satisfying~\eqref{eq:Sv}, and thus it halts and outputs that a.a.s.~$H_p$ contains a perfect matching.
We now show the correctness of this output.
We first show that \(M_0 \cup M_1\) can a.a.s.\ be extended to a matching \(M\) with \(\vec{i}_{\mathcal{P}}(V(H)\setminus V(M)) \in L^{s,\mu}_{\mathcal{P}}(H)\). 
For each \(\vec v\in I'\) and \(x\in S_{\vec v}\), choose a fixed collection of \(\eta n^{k-1}\) edges of \(H\) containing \(x\) and having index vector \(\vec v\), and expose only these edges. 
By the argument in Claim~\ref{claim:vec_v}, a.a.s.\ for every such pair \((\vec v,x)\), at least \((|V_0|+q)k\) of the exposed edges survive in \(H_p\) and are pairwise disjoint outside \(x\).  
We construct a matching \(M_3\) greedily by assigning to each \(v \in \bigcup_{\vec{v}\in I'} S_{\vec{v}}\) one such edge, avoiding previously chosen edges and \(M_0 \cup M_1\). 
This is possible since \(\sum_{\vec{v}} |S_{\vec{v}}|+\sum_{\vec{v}} c_{\vec{v}}^*\le q\), while each \(v\) has at least \((|V_0|+q)k\) available choices.
Thus we obtain a matching \(M_3\), disjoint from \(M_0 \cup M_1\), with \(\vec{i}_{\mathcal{P}}(V(M_3)) = \sum_{\vec{v}\in I'} (c_{\vec{v}} - c_{\vec{v}}^*) \vec{v}.\) 
Hence \(M:=M_0 \cup M_1 \cup M_3\) is a \((|V_0|+q)\)-solution.

\medskip
\textbf{From a solution to a perfect matching.}
Fix the matching \(M:=M_0\cup M_1\cup M_3\) constructed above, which is a \((|V_0|+q)\)-solution covering \(V_0\). 
It remains to show that a.a.s.~the remaining $k$-graph has a perfect matching. 
To this end, we verify that~\ref{item:(C*_degree)}--\ref{item:star} of Theorem~\ref{thm:auxiliary_thm_perfect_matching} hold for the remaining $k$-graph.

Let \(E'\) denote the set of all edges revealed so far. 
For every vertex \(v\in V(H)\setminus V_0\), at most \(\binom{d+k-1}{k}\eta n^{k-1}+|V_0|n^{k-2}+2\eta n^{k-1}\le \eta^{2/3}n^{k-1}\) edges of \(E'\) contain \(v\). 
We work with the \(k\)-graph obtained from \(H-V(M)\) by deleting all edges in \(E'\), and denote this graph by \(H'\). 
Let \(\mathcal P'=\{V'_0,V'_1,\dots,V'_d\}\) be the partition induced by \(\mathcal P\) on \(V(H')\), where \(V'_0=\emptyset\). 
All edges of \(H'\) are still unexposed, and hence the remaining random subgraph is distributed as \(H'_p\). 

We first verify the first assumption. Observe that for every vertex \(v \in V(H')\), at most 
\[
\frac{\binom{k}{\ell-1}\eta^{2/3}n^{k-1}}{\frac{\gamma }{2}\binom{|V(H')|-\ell}{k-\ell}}\le \sqrt{\eta}\binom{|V(H')|-1}{\ell-1}
\]
\((\ell-1)\)-sets \(S\) satisfy \(\deg_{H'}(S \cup \{v\}) < (\delta + \gamma/2)\binom{|V(H')|-\ell}{k-\ell}\) as \(\eta \ll \gamma\). 
This verifies~\ref{item:(C*_degree)}. 
Next, we verify the second assumption. 
Properties (1)--(2) of Lemma~\ref{lem:partition_l_degree} hold trivially for $\mathcal{P}'$, and we have $V_0'=\emptyset$. 
Moreover, for every \(i\in [s]\), $|V_i'|\ge \binom{c+k-2}{k-1}m$, and for every \(i\in \{s+1,\dots,d\}\), $|V_i'|\ge \delta' |V(H')|/4$. 
Recall that $H$ loses at most $k(|V_0|+q)$ vertices and for each vertex in $H'$, it loses at most $\eta^{2/3}n^{k-1}$ edges incident to it.
Thus, for every $\vec{v}\in I_{\mathcal{P},1}^{s,\mu}(H)$, there exist at least $\mu |V(H')|^k/2$ edges in $H'$ with index vector $\vec{v}$, and the same holds for every $\vec{v}\in I_{\mathcal{P},2}^{s,\mu}(H)$. 
It follows that \(I_{\mathcal{P}}^{s,\mu}(H)\subseteq I_{\mathcal{P}'}^{s,\mu/2}(H').\) 
Consequently, property (4) of Lemma~\ref{lem:partition_l_degree} is satisfied for $\mathcal{P}'$ (with slightly weakened parameters). 
It is easy to show that the last \(d-s\) parts of $\mathcal{P}'$ are $(\beta/2, 2^c)$-closed by \(\eta \ll \beta\), giving~\ref{item:(C*_partition)}. 
We now verify~\ref{item:star}. 
Let \(1\le i\le s\) and let
\(\vec{v}\notin I_{\mathcal P',2}^{s,\mu/2}(H')\)
satisfy \(\vec{v}|_i=1\) and
\(\sum_{j=s+1}^d \vec{v}|_j=k-1\).
Since \(I_{\mathcal P}^{s,\mu}(H)\subseteq
I_{\mathcal P'}^{s,\mu/2}(H')\), we also have
\(\vec{v}\notin I_{\mathcal P,2}^{s,\mu}(H)\). 
Therefore, by~\ref{item:extra_propert},
every vertex \(v\in V_i'\subseteq V_i\) is contained in fewer than \(\mu n^{k-1}\le 2\mu |V(H')|^{k-1}\) edges \(e\) with \(\vec{i}_{\mathcal P}(e)=\vec{v}\).

We verify the remaining assumption. 
Note that \(\vec{i}_{ \mathcal{P}'}(V(H'))=\vec{i}_{ \mathcal{P}}(V(H)\setminus V(M))\) and \(\vec{i}_{ \mathcal{P}}(V(H)\setminus V(M)) \in L_{ \mathcal{P}}^{s,\mu}(H)\), then \(\vec{i}_{ \mathcal{P}'}(V(H'))\in L_{ \mathcal{P}}^{s,\mu}(H)\subseteq L_{ \mathcal{P}'}^{s,\mu/2}(H').\)
Therefore, we can apply Theorem~\ref{thm:auxiliary_thm_perfect_matching} with parameters $\delta'/2$, $\mu/2$, $\sqrt{\eta}$, and $\beta/2$ in place of $\delta'$, $\mu$, $\eta$, and $\beta$, respectively, to $H'$.
We then conclude that a.a.s.~$H_p'$ contains a perfect matching.

\medskip
 \textbf{Running time.}
In the preprocessing step, Lemma~\ref{lem:partition_l_degree} finds the partition $\mathcal{P}$ in time $O(n^{2^ck+1})$. 
Computing all sets $V_{\vec{v}}^+$ and exposing edges to obtain $E(\vec{v})$ for all relevant $\vec{v}$ takes time $O(n^{k})$. 

In the search for a $(|V_0|+q)$-solution, we first enumerate all minimal matchings $M_0$, which takes time $O(n^{(k-1)|V_0|})$. 
Next, there are only constantly many choices for the coefficient vectors. 
For each $M_0$ and each such choice, we extend $M_0$ to a solution as follows: 
(i) selecting disjoint edge sets $M_{\vec{v}} \subseteq E(\vec{v})$, which takes time $O(n^{kq})$, and 
(ii) selecting pairwise disjoint vertex sets $S_{\vec{v}}$, which takes time $O(n^q)$.
Therefore, the total running time is $O(n^{(k-1)|V_0|+(k+1)q})$.

Note that $|V_0|$ is bounded from above by a constant depending only on $k$.
Therefore, there exists a constant $c=c(k)$ such that the overall running time is $O(n^{c})$.
\end{proof}

\subsection{Proof of Theorem~\ref{thm:auxiliary_thm_perfect_matching}}
We conclude this section by proving Theorem~\ref{thm:auxiliary_thm_perfect_matching}. 
To this end, we introduce the following useful constant.   
Given a set \( I \) of \( k \)-vectors in \( \mathbb{Z}^d \), and \( m \in \mathbb{N} \), consider the set \( J \) of all \( m' \)-vectors that are in the lattice generated by \( I \) with \( 0 \leq m' \leq m \). That is, for any \( \vec{v} \in J \), there exist \( a_{\vec{i}} \in \mathbb{Z}, \vec{i} \in I \) such that \(\vec{v} = \sum_{\vec{i} \in I} a_{\vec{i}} \vec{i}.\) 
Then \[ M^*(d, k, I, m):=\max_{\vec{v}\in J}\min\{\max_{\vec{i}\in I}|a_{\vec{i}}|:\vec{v} = \sum_{\vec{i} \in I} a_{\vec{i}} \vec{i}\},\] 
is the smallest integer \(M\) such that every \(\vec v\in J\) admits a
representation with \(|a_{\vec i}|\le M\) for all \(\vec i\in I\). 
Furthermore, let \(M(k, m)\) be the maximum of  \( M^*(d, k, I, m)\) over all \(d \leq d_0(k) := 2^{\binom{2k-1}{k-1}} + k\) and all families of \(k\text{-vectors } I \subseteq \mathbb{Z}^d.\)


We start with an outline of the proof.
Even though the partition \(\mathcal P=\{V_0,V_1,\dots,V_s,V_{s+1},\dots,V_d\}\) satisfies Lemma~\ref{lem:partition_l_degree} (1)-(5), it does not allow a na\"ive application of the random clustering argument in Theorem~\ref{thm:main_technical_theorem}, because the parts
\(V_1,\dots,V_s\) might be too small. 
To overcome this difficulty, we seek to first cover all vertices in \(V_1,\dots,V_s\) by a (random) matching, and then apply Theorem~\ref{thm:main_technical_theorem} to the remaining hypergraph. 

Note that such an application requires that the index vector of the remaining vertex set lies in \(L_{\mathcal P,1}^{s,\mu}(H)\), rather than just \(L_{\mathcal P}^{s,\mu}(H)\). 
This is where the simple lattice fact recorded as Fact~\ref{fact:separating_sublattice} in Section~\ref{sec:lattice_preparation} enters: we construct the first-stage matching with carefully prescribed index vectors so that the index vector of the remaining vertex set lies in \(L_{\mathcal P,1}^{s,\mu}(H)\), rather than merely in \(L_{\mathcal P}^{s,\mu}(H)\), before applying Theorem~\ref{thm:main_technical_theorem}.
To carry this out, we randomly partition \(U=\bigcup_{i=s+1}^d V_i\) into three parts \(U_1,U_2,U_3\). 
The first two parts are used to cover the vertices in \(V_1,\dots,V_s\). 
We first greedily cover almost all of these vertices using edges whose index vectors lie in \(I_{\mathcal P,2}^{s,\mu}(H)\). 
However, even after all vertices in \(V_1,\dots,V_s\) are covered, the index vector of the remaining vertex set need not lie in \(L_{\mathcal P,1}^{s,\mu}(H)\). 
Thus, for each \(i\in[s]\), we leave a bounded number of vertices in \(V_i\) uncovered in the greedy step, and cover them later using carefully prescribed index vectors. 
These prescriptions are obtained by expressing the current residue as an integer combination of vectors from \(I_{\mathcal P,1}^{s,\mu}(H)\) and \(I_{\mathcal P,2}^{s,\mu}(H)\), and then choosing the last few edges so that the \(L_{\mathcal P,2}^{s,\mu}(H)\)-component is eliminated. 
The third part \(U_3\) is kept unused during this correction step. 
This ensures that, after the vertices in \(V_1,\dots,V_s\) have been covered, the remaining hypergraph still satisfies the degree, closedness and lattice assumptions needed for Theorem~\ref{thm:main_technical_theorem}. 
Applying that theorem gives a spread distribution on perfect matchings of the remaining hypergraph. 
The spread property of the final matching follows by combining this distribution with the random choices made in the first stage: each edge used to cover a vertex in \(V_1,\dots,V_s\) is chosen uniformly from many available choices, and hence is selected with probability \(O(1/n^{k-1})\).

\begin{proof}[Proof of Theorem~\ref{thm:auxiliary_thm_perfect_matching}]
Note that \(c_{k,\ell}^*\ge c_{k,k-1}^*=1/k \), then \(\lfloor 1/c_{k,\ell}^* \rfloor\le k\). It follows that \(d-s\le k\), \(|I_{\mathcal{P},2}^{s,\mu}(H)|\le s\binom{2k-2}{k-1}\), and \(s\le 2^{\binom{2k-2}{k-1}}\). Let \(q:=(2k+1)^k\). 
We choose \[t:=2^k\qquad\text{and}\qquad m=M(k,qk^2).\]
Define additional constants \(\alpha,\varepsilon>0\) such that
\[
1/n_0 \ll \eta\ll 1/C''\ll\mu,\beta\ll\alpha,\varepsilon \ll \delta' \ll \gamma, c_{k,\ell}^*,1/k.
\] 
Suppose that \(H\) and \(\mathcal{P}\) satisfy the assumptions of the theorem. 
Let \(U=\bigcup_{ i\in \{s+1,\dots,d\}}V_i\). 
We first randomly partition \(U\) into three parts of almost equal size, denoted by \(U_1, U_2,\) and \(U_3\). 
By applying Lemma~\ref{lem:application_McDiarmid} to each of the following properties and taking a union bound, with probability at least \(1-\exp(-\Omega(n))\), all of the following hold:
\begin{enumerate}[label=(U\arabic*)]
    \item For every \(v\in \bigcup_{i\in [s]}V_i\), we have \(\deg_H(v,U_1)\ge \gamma n^{k-1}\);\label{item:U_1}
    \item For every \(i \in [s]\) and $\vec{v} \in I_{\mathcal{P},2}^{s,\mu}(H)$ with $\vec{v}|_{i} = 1$, every vertex \(v\in V_i\) is in at least \(\mu^2n^{k-1}\) edges \(e\) with index vector  \(\vec{v}\) such that all vertices of \(e\) except \(v\) lie in \(U_2\);\label{item:U_2}
    \item For every \(i\in \{s+1,\dots,d\}\), we have \(|V_i\cap U_3|\ge \delta'n/10\), and for every \(u,v\in V_i\), there are at least \(\beta^2 n^{tk-1}\) reachable sets \(S\subseteq U_3\) for \(u\) and \(v\);\label{item:U_3_reachable}
    \item For every $\vec{u} \in I_{\mathcal{P},1}^{s,\mu}(H)$, there are at least \(\mu^2n^{k}\) edges in \(H[U_3]\) with index vector \(\vec{u}\).\label{item:U_3_lattice}
\end{enumerate}
Fix a partition satisfying the above properties. 

\medskip
\textbf{Step 1: Covering all vertices in \(\bigcup_{i\in [s]}V_i\).} 
First, we construct a matching \(\tilde{M}_1\) using only vertices in \(U_1\) to cover almost all vertices in \(\bigcup_{i\in [s]}V_i\). 
For every \(i\in [s]\), let \(v_1,\dots ,v_{|V_i|}\) be an enumeration of \(V_i\), and let \(I_i\) be the ordered collection of vectors \(\vec{v}\) in \(I_{\mathcal{P},2}^{s,\mu}(H)\) with \(\vec{v}|_i=1\). 
We cover the vertices \(v_1,\dots,v_{|V_i|-m|I_i|}\) using a random greedy algorithm: process these vertices in order, and for each vertex \(v_i\), choose uniformly at random an edge containing \(v_i\) and \(k-1\) vertices in \(U_1\), such that the edge has index vector in \(I_{\mathcal{P},2}^{s,\mu}(H)\) and all chosen edges are pairwise disjoint. 
Now note that for $i\in [s]$, each $V_i$ has $m|I_i|$ vertices not covered by $\tilde M_1$.

Next, we extend \(\tilde{M}_1\) to a matching \(M_1\) that covers the remaining vertices in \(\bigcup_{i\in[s]}V_i\), such that the other \(k-1\) vertices of each added edge lie in \(U_2\). 
Set \(\vec{x}:=\vec{i}_{\mathcal{P}}\bigl(V(H)\setminus V(\tilde{M}_1)\bigr).\)
Since \(\vec{i}_{\mathcal P}(V(H))\in L_{\mathcal P}^{s,\mu}(H)\) and every edge of \(\tilde M_1\) has index vector in \(I_{\mathcal P,2}^{s,\mu}(H)\), we have \[\vec{x}\in L_{\mathcal{P}}^{s,\mu}(H)=L_{\mathcal{P},1}^{s,\mu}(H)+ L_{\mathcal{P},2}^{s,\mu}(H).\] 
By Fact~\ref{fact:separating_sublattice}, there exists \(\vec{y}\in L_{\mathcal{P},2}^{s,\mu}(H)\) such that \(\vec{x}-\vec{y}\in L_{\mathcal{P},1}^{s,\mu}(H)\).
To apply Theorem~\ref{thm:main_technical_theorem}, we need to choose such a vector \(\vec{y}\) that can be realized as the total index vector of the edges added to \(\tilde{M}_1\). 
We now construct such a vector. 

Recall that, for every \(i\in [s]\), the \(i\)-th coordinate of \(\vec{x}\) is exactly \(m|I_i|\). 
As each $\vec{v}\in I_i$ satisfies $\vec{v}\in I_{\mathcal{P},2}^{s,\mu}(H)$, we get
\[\vec{x}-\sum_{i\in [s]}\sum_{\vec{v}\in I_i} m\vec{v}\in L_{\mathcal{P}}^{s,\mu}(H).\]
Moreover, for each \(i\in \{s+1,\dots,d\}\), by the pigeonhole principle, there exists an integer \(q_i \le q\) such that 
\[
\vec{u}_i := (0,\dots,0, q_i k, 0,\dots,0) \in L_{\mathcal{P},1}^{s,\mu}(H),
\] where the \(i\)-th coordinate is \(q_i k\) and all other coordinates are zero. 
Indeed, otherwise these $q+1$ multiples of the vector $(0,\dots,0, k, 0,\dots,0)$ would lie in distinct cosets of \(L_{\mathcal P,1}^{s,\mu}(H)\) in \(L_{\max}^{d-s}\), contradicting Lemma~\ref{lem:coset_group_of_l_degree}.
Choosing appropriate integers \(c_{\vec{u}_i}\), we obtain
\begin{equation}\label{eq:w}
    \vec{w}:= \vec{x}-\sum_{i\in [s]}\sum_{\vec{v}\in I_i} m\vec{v}- \sum_{i\in [d]\setminus [s]} c_{\vec{u}_i}\vec{u}_i \in L_{\mathcal{P}}^{s,\mu}(H),
\end{equation}
such that \(\vec{w}|_i =0\) for all \(i\in [s]\) and \(0\le \vec{w}|_i \le q_i k\) for all \(i\in \{s+1,\dots,d\}\). 
Let \(I^{(1)} := I_{\mathcal{P},1}^{s,\mu}(H)\) and \(I^{(2)} := I_{\mathcal{P},2}^{s,\mu}(H)\). 
Since the sum of the coordinates of \(\vec{w}\) is at most \(qk(d-s)\le qk^2\), by the definition of \(m=M(k,qk^2)\), we may write
\[\vec{w}=\sum_{\vec{u}\in I^{(1)}}a_{\vec{u}}\vec{u} +\sum_{\vec{v}\in I^{(2)}} b_{\vec{v}}\vec{v},\]
where \(a_{\vec{u}},b_{\vec{v}}\in\mathbb{Z}\) and \(|b_{\vec{v}}|\le m\) for every \(\vec{v}\in I^{(2)}\).
Since \(b_{\vec{v}}\) may be negative, write \(b_{\vec{v}}=b_{\vec{v}}^+-b_{\vec{v}}^-\), where \(b_{\vec{v}}^+,b_{\vec{v}}^-\) are non-negative integers and at least one of them is zero. 
Grouping the vectors in \(I^{(2)}\) according to the collections \(I_i\), we have
\[ 
\vec{w}=\sum_{\vec{u}\in I^{(1)}}a_{\vec{u}}\vec{u}+ \sum_{i\in [s]}\sum_{\vec{v}\in I_i}(b^+_{\vec{v}}-b^-_{\vec{v}})\vec{v}.
\]
We now define
\[\vec{y}:=\sum_{i\in[s]}\sum_{\vec{v}\in I_i}(m+b_{\vec{v}}^+-b_{\vec{v}}^-)\vec{v}.\]
Clearly, \(\vec{y}\in L_{\mathcal{P},2}^{s,\mu}(H)\). 
Combining~\eqref{eq:w} with the above decomposition of \(\vec{w}\), we obtain
\[
\vec{x}-\vec{y}= \sum_{\vec{u}\in I^{(1)}}a_{\vec{u}}\vec{u}+\sum_{i\in [d]\setminus [s]} c_{\vec{u}_i}\vec{u}_i
\in L_{\mathcal{P},1}^{s,\mu}(H).
\]
Thus, \(\vec{y}\) is a correction vector of the type described in Fact~\ref{fact:separating_sublattice}, and it remains to realize \(\vec y\) as the total index vector of a matching. 

The above decomposition determines the required number of edges of each index vector. 
Since \(|b_{\vec{v}}|\le m\), we have \(m+b_{\vec{v}}^+-b_{\vec{v}}^-=m+b_{\vec{v}}\ge0\). 
Thus, previously chosen edges do not need to be modified, and the remaining lattice discrepancy can be corrected by selecting additional edges of the appropriate types. 
In particular, for every \(\vec{v}\in I_i\), we need \(m+b_{\vec{v}}^+-b_{\vec{v}}^-\) edges with index vector \(\vec{v}\). 
Note that these quantities sum to \(m|I_i|\), so this exactly matches the number of uncovered vertices in \(V_i\). 
We now turn this requirement into an explicit construction. 

We proceed cluster by cluster. 
For each \(i\in [s]\), 
we construct the edges for \(V_i\) by processing the vectors in \(I_i\) sequentially. 
For each \(\vec{v}\in I_i\), we select \(m+b^+_{\vec{v}}-b^-_{\vec{v}}\) vertices from the remaining uncovered vertices of \(V_i\), and for each such vertex \(v\), choose uniformly at random an edge containing \(v\) and \(k-1\) vertices in \(U_2\), whose index vector is \(\vec{v}\), ensuring that all chosen edges are pairwise disjoint.
We then proceed similarly for all remaining vectors in $I_i$, $i\in [s]$. 
This yields a matching \(M_1\) covering all vertices in \(\bigcup_{i\in [s]} V_i\) such that \[\vec{i}_{\mathcal{P}}(V(H)\setminus V(M_1))=\vec{x}-\vec{y} \in L_{\mathcal{P},1}^{s,\mu}(H).\]

This whole procedure can be carried out.
Indeed, for the first phase, for every vertex in \(\bigcup_{i\in [s]}V_i\), the number of edges containing it whose index vector does not lie in \(I^{(2)}\) is at most \(s\binom{2k-2}{k-1}\cdot4\mu n^{k-1}\) by assumption \ref{item:star}.  
Therefore, for each of the first \(|V_i|-m|I_i|\) vertices, the number of available choices is at least
\[\gamma n^{k-1}-s\binom{2k-2}{k-1}\cdot 4\mu n^{k-1}-k\delta'n\cdot n^{k-2}\ge \gamma n^{k-1}/2,\]
where we use that the partition satisfies \ref{item:U_1} and that
\(1/n\ll \mu \ll \delta' \ll \gamma,1/k\). 
Moreover, for the second phase, that is, for the last \(m|I_i|\) vertices in each cluster, \ref{item:U_2} implies that the number of available choices is at least \[\mu^2n^{k-1}-ms\binom{2k-2}{k-1}k\cdot n^{k-2}\ge \mu^2 n^{k-1}/2\] 
using that \(\sum_{i\in [s]}|I_i|=|I^{(2)}|\le s\binom{2k-2}{k-1}\). 

\medskip
\textbf{Step 2: Finding a spread matching in the remaining hypergraph.} 
Since all vertices in the first \(s\) parts of \(\mathcal{P}\) are covered by \(M_1\), we restrict attention to the remaining parts. 
Next, consider the remaining hypergraph \(H':=H[V(H)\setminus V(M_1)]\) and its partition \(\mathcal{P}':=\{V_i\cap V(H')\}_{i\in \{s+1,\dots,d\}}\). 
Set \(n':=|V(H')|\). 
Note that \(n'\ge n-k\cdot c\delta'n>n/2\). 
We verify that \(H'\) and \(\mathcal{P}'\) satisfy the assumptions of Theorem \ref{thm:main_technical_theorem}. 

We first verify assumption~\ref{item:almost_perfect_degree} of Theorem~\ref{thm:main_technical_theorem}.
Recall that $\delta(k,\ell,2k-\ell-1)=c_{k,\ell}^*$ by Theorem \ref{thm:deltaFD}. 
Fix a vertex \(v\in V(H')\). 
By assumption~\ref{item:(C*_degree)}, all but at most \(\eta\binom{n-1}{\ell-1}\) sets \(S\in\binom{V(H)\setminus\{v\}}{\ell-1}\) satisfy \(\deg_H(S\cup\{v\})\ge (c_{k,\ell}^*+\gamma)\binom{n-\ell}{k-\ell}.\) 
Now let \(S\subseteq V(H)\setminus\{v\}\) be such a set. Note that \[\deg_{H'}(S\cup \{v\})\ge (c_{k,\ell}^*+\gamma)\binom{n-\ell}{k-\ell}-k\cdot c\delta'n \cdot n^{k-\ell-1}\ge (c_{k,\ell}^*+\gamma/2)\binom{n'-\ell}{k-\ell}.\] 
Therefore, all but at most \(\eta\binom{n-1}{\ell-1}\le 2\eta\binom{n'-1}{\ell-1}\) sets \(S\in\binom{V(H')\setminus\{v\}}{\ell-1}\) satisfy \(\deg_{H'}(S\cup \{v\})\ge (c_{k,\ell}^*+\gamma/2)\binom{n'-\ell}{k-\ell}\), as required. 
Next, we verify the assumption~\ref{item:mathcal_P} of Theorem~\ref{thm:main_technical_theorem}. 
By \ref{item:U_3_reachable} and the fact that \(U_3 \subseteq V(H')\), it follows that 
\(\mathcal{P}'\) is a \(\bigl(\beta^2,\, t,\, \delta'/10\bigr)\)-good partition of \(V(H')\).
We now verify the assumption~\ref{item:almost_robust_edge} of Theorem~\ref{thm:main_technical_theorem}. 
Note that 
\[
\delta_1(H') \ge \binom{k-1}{\ell-1}^{-1} {(1-2\eta)\binom{n'-1}{\ell-1} \cdot (c_{k,\ell}^*+\gamma/2)\binom{n'-\ell}{k-\ell}}\ge 2\gamma {(n')}^{k-1}.
\]
Let $ B_0 $ be the set of all vertices $v\in V(H')$ such that \(\card{E_{\mathcal{P}'}^{\mu^2}(v)} \le \varepsilon {(n')}^{k-1}\).
Note that for each \(v\in B_0\), the vertex \(v\) lies in at least \(2\gamma {(n')}^{k-1} - \varepsilon {(n')}^{k-1} \ge \gamma {(n')}^{k-1}\) edges whose index vector is not in \(I_{\mathcal{P}'}^{\mu^2}(H')\), yielding
\[
|B_0| \le {k\binom{2k-1}{k}\mu^2 {(n')}^k}/{\gamma {(n')}^{k-1}}<\alpha {n'},
\]
as $\mu\ll\alpha,\gamma, 1/k$. 
Finally, we verify the fourth assumption. 
By \ref{item:U_3_lattice}, we have \(L_{\mathcal{P},1}^{s,\mu}(H)\subseteq L_{\mathcal{P}'}^{\mu^2}(H')\). 
Since \(|\mathcal{P}'| = d-s\), by Lemma~\ref{lem:coset_group_of_l_degree}, it follows that 
\[|Q(\mathcal{P}',L_{\mathcal{P}'}^{\mu^2}(H'))|=|L_{\max}^{d-s}/L_{\mathcal{P}'}^{\mu^2}(H')|\le |L_{\max}^{d-s}/L_{\mathcal{P},1}^{s,\mu}(H)|\le q.\] 
Moreover, we have \(\vec{i}_{\mathcal{P}'}(V(H'))=\vec{i}_{\mathcal{P}}(V(H)\setminus V(M_1)) \in L_{\mathcal{P},1}^{s,\mu}(H)\subseteq L_{\mathcal{P}'}^{\mu^2}(H')\).
Therefore, conditioning on any instance of \(M_1\), we apply Theorem \ref{thm:main_technical_theorem} to \(H'\) and \(\mathcal{P}'\) to find an \(\left(C''/(n')^{k-1}\right)\)-spread random perfect matching \(M_2\) in \(H'\). 
Set \(M=M_1\cup M_2\), which defines a probability distribution on the perfect matchings of $H$.

It remains to show that this distribution is \(O(1/n^{k-1})\)-spread.  
Let \(S \subseteq E(H)\) be a matching and we show that \(\mathbb{P}[S \subseteq M] = (O(1/n^{k-1}))^{|S|}\).  
Note that if \(S \subseteq M\), then every edge in \(S\) must either contain exactly one vertex in \(\bigcup_{i\in [s]} V_i\), or be entirely contained in \(U\); otherwise, \(\mathbb{P}[S \subseteq M]=0\). 
Let \(S_1\) be the set of edges in \(S\) that contain exactly one vertex in \(\bigcup_{i\in [s]} V_i\), and let \(S_2 := S \setminus S_1\). Then
\[\mathbb{P}[S\subseteq M ]=\mathbb{P}[S_1\subseteq M_1]\cdot \mathbb{P}[S_2 \subseteq M_2 \mid S_1 \subseteq M_1].\]
On the one hand, \(S_1\subseteq M_1\) only if for every edge \(e \in S_1\), the edge chosen to cover the unique vertex in 
\(e \cap \left(\bigcup_{i\in [s]} V_i\right)\) is exactly \(e\). 
Since each such choice was made uniformly at random from at least \(\mu ^2n^{k-1}/2\) possibilities, it follows that \(\mathbb{P}[S_1\subseteq M_1 ]=(O(1/n^{k-1}))^{|S_1|}.\) 
On the other hand, by construction, \(M_2\) is \(O(1/(n')^{k-1})\)-spread conditioned on any instance of \(M_1\), and hence \(\mathbb{P}[S_2 \subseteq M_2 \mid S_1 \subseteq M_1]
= (O(1/n^{k-1}))^{|S_2|}.\) 
Therefore, \(\mathbb{P}[S \subseteq M] = (O(1/n^{k-1}))^{|S|}\), as desired. 
The final assertion follows from Theorem~\ref{thm:FKNP}.
\end{proof}

\section{The enumeration aspect}\label{sec:enumeration}

    Corollary~\ref{lem:counting_PM} and Corollary~\ref{lem:counting_packing} follow from the fact below together with Theorem~\ref{thm:auxiliary_thm_perfect_matching} and Theorem~\ref{thm:main_technical_theorem}, respectively. 
    More precisely, if \(H\) satisfies the corresponding divisibility conditions (as in Theorem~\ref{thm:gan_han_structural} or Theorem~\ref{thm:han_structural_theorem}), then the corresponding theorem yields a spread distribution in the subgraph \(H' := H - V(M_0)\), where \(M_0\) is the associated solution.  
    Define \(\mathcal{M}(H)\) to be the set of perfect \(F\)-packings of \(H\). 
    By the fact below, it follows that \(
    \card{\mathcal{M}(H')}\ge\left(\varepsilon (n-c)\right)^{\frac{e_F}{r m_1(F)}(n-c)},\) where \(c=O(1)\) is the number of vertices covered by \(M_0\).  
    Observe that every perfect \(F\)-packing \(M\) in \(H'\) can be extended to a perfect \(F\)-packing \(M \cup M_0\) in \(H\). 
    Hence, \(\card{\mathcal{M}(H)}\ge\card{\mathcal{M}(H')}\ge (\varepsilon n/2)^{\frac{e_F}{r m_1(F)}n}\).
    
    \begin{fact}
        Let \(F\) be an \(r\)-vertex \(k\)-graph, and \(H\) an \(n\)-vertex \(k\)-graph. If there is a \(\left(C/n^{1/m_1(F)}\right)\)-spread distribution \(\textbf{M}\) on the set of $F$-factors in \(H\), then $H$ contains at least $(\varepsilon n)^{\frac{e_F}{rm_1(F)}n}$ $F$-factors for \(\varepsilon = C^{-m_1(F)}\).
    \end{fact}
    \begin{proof}
        Fix an \(F\)-factor \(M\in \mathcal{M}(H)\), which contains \(e_Fn/r\) edges as a subgraph of \(H\). Since \(\textbf{M}\) is \(\left(C/n^{1/m_1(F)}\right)\)-spread, we have \(\mathbb{P}[\textbf{M}=M]\le (C/n^{1/m_1(F)})^{e_Fn/r}\). 
        Thus, \[1=\sum_{M\in \mathcal{M}(H)}\mathbb{P}[\textbf{M}=M]\le \card{\mathcal{M}(H)}(C/n^{1/m_1(F)})^{e_Fn/r},\]
        which implies \(\card{\mathcal{M}(H)}\ge (n^{1/m_1(F)}/C)^{e_Fn/r}= (\varepsilon n)^{\frac{e_F}{rm_1(F)}n}\), as desired. 
    \end{proof}

\section{Concluding remarks}\label{remark}

In this paper, we established an algorithmic result (Theorem~\ref{thm:algorithm_PM}) that guarantees, with high probability, the existence of a perfect matching in a random subgraph of \(H\) when the algorithm accepts, and the non-existence of such a matching when it rejects. 
 

In contrast to the perfect matching setting, we are unable to obtain an analogous algorithm for \(F\)-packings as in Theorem~\ref{thm:algorithm_PM}. 
Indeed, following the proof strategy of Theorem~\ref{thm:algorithm_PM}, it may be necessary to expose many copies of \(F\) in the random subgraph in order to verify certain divisibility conditions. 
Even if, for each vertex \(v\), we only expose at most \(\eta n^{r-1}\) copies of \(F\) containing \(v\), each such copy may involve several edges containing \(v\), so the exposure process may reveal almost all edges incident to a vertex as non-edges.
This prevents us from applying Theorem~\ref{thm:main_technical_theorem} to extend a suitable solution to a perfect \(F\)-packing.
In the perfect matching setting, by contrast, \(F\) itself consists of a single edge, so the exposure process only reveals individual edges, and the loss of degree can be quantitatively controlled throughout the algorithm. 

We note that while our main approach follows the randomized algorithm of Kelly, M\"{u}yesser, and Pokrovskiy, we also explored alternative methods based on iterative absorption, as developed by Pham, Sah, Sawhney, and Simkin~\cite{pham2022toolkit}. 
While we think both approaches work, we chose to present the proof using the method of Kelly et al.~as it is conceptually simpler.

\section*{Acknowledgment}
We would like to thank Yu Chen and Pan Peng for helpful discussions on this project.

\bibliographystyle{abbrv}
\bibliography{main}

\end{document}